\documentclass[11pt,twoside]{amsart}

\usepackage{amscd}
\usepackage[pdftex]{hyperref}
\hypersetup{colorlinks,citecolor=blue,linktocpage,hyperindex=true,backref=true}

\newtheorem{thm}{Theorem}[section]
\newtheorem{lem}[thm]{Lemma}
\newtheorem{cor}[thm]{Corollary}

\theoremstyle{definition}

\newtheorem{defn}[thm]{Definition}
\newtheorem{say}[thm]{}
\newtheorem{exmp}[thm]{Example}

\newtheorem{rem}[thm]{Remark}	

\theoremstyle{remark}

\newtheorem{claim}[thm]{Claim}
\newtheorem{case}{Case}
\newtheorem{subcase}{Subcase}
\newtheorem{step}{Step}

\newcommand{\kbys}{K+B^{Y^{(s)}}}
\newcommand{\kbxs}{K+B^{X^{(s)}}}

\newcommand{\supp}{\operatorname{supp}}

\newcommand{\xso}{x^{(s)}_{1}}
\newcommand{\xst}{x^{(s)}_{2}}
\newcommand{\xsk}{x^{(s)}_{k}}

\newcommand{\xstk}{x^{(s)}_{2;k}}
\newcommand{\xstkz}{x^{(s)}_{2;k_{0}}}

\newcommand{\xspk}{x^{(s)}_{p;k}}

\newcommand{\rk}{\operatorname{rk}}
\newcommand{\Pic}{\operatorname{Pic}}
\newcommand{\Aut}{\operatorname{Aut}}
\newcommand{\Exc}{\operatorname{Exc}}
\newcommand{\Char}{\operatorname{char}}

\title
[Boundedness and
$K^{2}$ for Log Surfaces]
{Boundedness and
$\bold K^{\boldsymbol 2}$ for Log Surfaces\\
}
\author{Valery Alexeev}
\date{April 1, 1994}

\begin{document}
\maketitle
\bibliographystyle{amsplain}

\tableofcontents

\setcounter{section}{-1}
\section{Introduction}
\label{Introduction}

  \begin{say}
The aim of this paper is to generalize to the singular case the
following statements which in the nonsingular case are easy and well-known:
  \begin{enumerate}
    \item The class of smooth
          surfaces with ample anticanonical divisor $-K$, also known as Del
          Pezzo surfaces, is bounded.
    \item The class of smooth
          surfaces with ample canonical divisor $K$ with $K^{2}\le
	  C$ is bounded.
    \item The class of smooth
          surfaces with ample canonical divisor $K$ with $K^{2}=
	  C$ is bounded.
    \item The set $\{K_{X}^{2}\}$, where $X$ is a smooth surface with
          ample $K_{X}$, is the set $\mathbb N$ of positive integer numbers.
  \end{enumerate}
  \end{say}

  \begin{say}
Of course, the last of these statements is trivial, and the second and the
third ones are equivalent. In the singular case, however, the
set $\{K_{X}^{2}\}$ {\em a priori\/} is just a certain subset in the
set of positive
rationals. We shall prove that under natural and rather weak
conditions this set has a rigid structure: it satisfies the descending
chain condition, abbreviated below to
D.C.C. for short. The boundedness for surfaces with
constant $K^{2}$ will be proved under much weaker conditions on
singularities than for surfaces with $K^{2}$ only bounded.
  \end{say}

  \begin{say}
In order to be useful, the generalizations should come from interesting
real-life examples. The examples should also suggest what conditions on
singularities are most natural.
The section~\ref{Applications} contains several applications, due to
J.Koll\'ar,
G.Xiao and others, that provide the motivation for our boundedness
theorems.
The most important application is the
projectiveness of the coarse moduli space of stable surfaces.
Two others are a formula bounding the automorphism group of a
(possibly, singular) surface of general type by $cK^{2}$, and
a theorem on the uniform plurigenera of elliptic threefolds.
  \end{say}

  \begin{say}
\label{proved in this paper}

These are the precise statements of the results proved in this paper:
  \begin{enumerate}
\label{proved here}
 \item Fix $\varepsilon>0$.
          Consider all projective surfaces $X$
          with an $\mathbb R$-divisor $B=\sum b_{j}B_{j}$ such that
          $K_{X}+B$ is MR $\varepsilon$-log terminal and $-(K_{X}+B)$ is
          nef excluding only those for which at the same time $K_{X}$
          is numerically trivial, $B$ is zero and $X$ has at worst Du Val
          singularities.
          Then the class $\{X\}$ is bounded
          (theorem~\ref{bound for general when -K nef}).
    \item Fix $\varepsilon>0$, a constant $C$ and a D.C.C. set $\mathcal C$.
          Consider all surfaces $X$
          with an $\mathbb R$-divisor $B=\sum b_{j}B_{j}$ such that
          $K_{X}+B$ is MR $\varepsilon$-log terminal, $K_{X}+B$ is
          big and nef, $b_{j}\in \mathcal C$ and $(K+B)^{2}\le C$.
          Then the class $\{(X,\supp B)\}$ is bounded
          (theorem~\ref{bound for varepsilon-terminal, K big, nef}).
    \item Fix a constant $C$ and a D.C.C. set $\mathcal C$.
          Consider all surfaces $X$
          with an $\mathbb R$-divisor $B=\sum b_{j}B_{j}$ such that
          $K_{X}+B$ is MR semi-log canonical, $K_{X}+B$ is
          ample, $b_{j}\in \mathcal C$ and $(K+B)^{2}= C$.
          Then the class $\{(X,\sum b_{j}B_{j})\}$ is bounded
          (theorem~\ref{K2=const for semi-log canonical}).
    \item Fix a D.C.C. set $\mathcal C$.
          Consider all surfaces $X$
          with an $\mathbb R$-divisor $B=\sum b_{j}B_{j}$ such that
          $K_{X}+B$ is MR semi-log canonical, $K_{X}+B$ is
          ample and $b_{j}\in \mathcal C$.
          Then the set $\{ (K_{X}+B)^{2} \}$ is a D.C.C. set
          (theorem~\ref{DCC for semi-log canonical}).
  \end{enumerate}
  \end{say}

  \begin{say}
The necessary definitions will be given in the section
{}~\ref{Standard definitions}. As the
first approximation, the reader can
look at the results obtained by dropping $B$. It is
interesting to note that none of
the conditions above can be weakened, see the examples in section
{}~\ref{Applications}.
(Of course, it is always possible to choose an entirely different way
of generalizing the nonsingular case and pick a set of conditions
which is orthogonal to ours.)
  \end{say}

  \begin{say}
I became interested in the subject reading
\cite{Kollar92} by Koll\'ar, and that is where most of the motivation
for the present paper comes from. The statements (3) and (4) above solve
and generalize a Koll\'ar's conjecture formulated there. A hope for
proving this conjecture came to me when I looked at the preprint
\cite{Xiao91} by Xiao, and especially at his proof of Proposition 5. I
am most indebted to G.Xiao for answering many of my questions
concerning his preprint. I am also thankful to V.V.Shokurov and
K.Matsuki for useful discussions.
  \end{say}

  \begin{say}
Throughout, we work with projective algebraic schemes over an
algebraically closed field of arbitrary
characteristic.
  \end{say}

\section{Standard definitions}
\label{Standard definitions}

  \begin{say}
For a normal variety $X$, $K_{X}$ or simply $K$ will always denote the
class of linear equivalence of the canonical Weil divisor.

  \end{say}

  \begin{defn}
An {\em $\mathbb R$-divisor\/} $D=\sum d_{i}D_{i}$
is a linear combination
of prime Weil divisors with real coefficients,
i.e.\ an element of $N^{1}\otimes \mathbb R$.
 An $\mathbb R$-divisor is said to be
$\mathbb R$-Cartier if it is a combination
of Cartier divisors with real coefficients, i.e.\ if it belongs to
$Div(X)\otimes\mathbb R \subset N^{1}(X)\otimes\mathbb R$.
The $\mathbb Q$-divisors and $\mathbb Q$-Cartier
divisors are defined in a similar fashion.
  \end{defn}

  \begin{defn}
Consider an $\mathbb R$-divisor
$K+B=K_{X}+\sum b_{j}B_{j}$ and assume that
    \begin{enumerate}
      \item $K+B$ is $\mathbb R$-Cartier
      \item $0\le b_{j}\le1$
    \end{enumerate}

For any resolution $f:Y\to X$ look at the natural formula
      \begin{eqnarray}
\label{definition of codiscrepancies}
K_{Y}+B^{Y}=
f^{*}(K_{X}+\sum b_{j}B_{j})= K_{Y}+\sum b_{j}f^{-1}B_{j} + \sum b_{i}F_{i}
      \end{eqnarray}
or, equivalently,

      \begin{eqnarray}
\label{definition of log discrepancies}
K_{Y}+\sum b_{j}f^{-1}B_{j} + \sum F_{i}=
f^{*}(K_{X}+\sum b_{j}B_{j}) + \sum f_{i}F_{i}
      \end{eqnarray}

Here $f^{-1}B_{j}$ are the proper
preimages of $B_{j}$ and $F_{i}$ are the exceptional divisors of
$f:Y\to X$.
The coefficients $b_{i},b_{j}$ are called codiscrepancies, the
coefficients $f_{i}=1-b_{i},f_{j}=1-b_{j}$ -- log discrepancies.

  \end{defn}

  \begin{rem}
In fact, $K+B$ is not a usual $\mathbb R$-divisor but rather a special
gadget consisting of a linear class of a Weil divisor $K$ (or a
corresponding reflexive sheaf) and an honest $\mathbb R$-divisor $B$.
This, however, does not cause any confusion.
  \end{rem}

  \begin{defn}
An $\mathbb R$-Cartier divisor $K+B$ (or a pair $(X,B)$) is called
    \begin{enumerate}
      \item log canonical, if the log discrepancies $f_{k}\ge0$
      \item Kawamata log terminal, if $f_{k}>0$
      \item canonical, if $f_{k}\ge0$
      \item terminal, if $f_{k}>0$
      \item $\varepsilon$-log canonical, if $f_{k}\ge\varepsilon$
      \item $\varepsilon$-log terminal, if $f_{k}>\varepsilon$
    \end{enumerate}
for every resolution $f:Y\to X$.
  \end{defn}

  \begin{say}
The name $\varepsilon$-log terminal was suggested to me by V.V.Shokurov.
  \end{say}

  \begin{defn}
In the two-dimensional case
we shall say that $K+B$ is MR log canonical, MR Kawamata log terminal
etc. if we require the previous inequalities to hold not for all
resolutions $f:Y\to X$ but only for a distinguished one,
the minimal desingularization.
  \end{defn}

  \begin{say}
In the surface case
all the standard theorems of the log Minimal Model Program are valid under
very weak assumptions, see the section~\ref{MMP in dimension 2}. In
particular, assuming $K+B$ to be MR Kawamata log terminal or MR log
canonical is sufficient for all applications.
  \end{say}

  \begin{defn}
Let $X$ be a reduced scheme satisfying the Serre's condition $S_{2}$ and
$B=\sum b_{j}B_{j}$, $0\le b_{j} \le 1$ be an $\mathbb R$-divisor.
An $\mathbb R$-Cartier divisor $K_{X}+B$  (or a pair $(X,B)$) is called
semi-log canonical if $f_{k}\ge0$ in the above formula
{}~\ref{definition of log discrepancies} for every {\em semiresolution}
$f:Y\to X$ (see \cite{FAAT} 12.2.1).
\end{defn}

\begin{say}
In dimension 2 and characteristic 0 the semi-log canonical
singularities are classified in \cite{KollarShepherdBarron88}. This is
the complete list (modulo analytic isomorphism): nonsingular points,
cones over nonsingular elliptic curves, double normal crossing points
$xy=0$, pinch points $x^{2}=y^{2}z$; and all finite quotients of above.
In codimension one the semi-log canonical schemes have only normal
crossing points.
\end{say}

\begin{say}
\label{reduce semi-log canonical to log canonical}

Let $\nu: X^{\nu}\to X$ be a normalization of $X$,
$X^{\nu}=\cup X_{m}$ be a decomposition into irreducible components
and define $B_{m}$ on $X_{m}$ to be $\nu^{-1}(B)$ plus the double
intersection locus. Then
  \begin{displaymath}
K_{X^{\nu}}+\sum B_{m}= \nu^{*}(K_{X}+B)
  \end{displaymath}
and $K_{X_{m}}+B_{m}$ are log canonical ( \cite{FAAT} 12.2.2,4).
\end{say}

\begin{defn}
We shall say that $K_{X}+B$ is MR semi-log canonical if $K_{X}$ is
semi-log canonical and all $K_{X_{m}}+B_{m}$ are MR log canonical.
\end{defn}

\section{Examples}
\label{Applications}

  \begin{say}
I include the following three examples
from \cite{Kollar92}, in a  generalized form.
They help to understand
what kind of boundedness results for singular surfaces are desirable.
After the theorems of this paper are proved, these applications are no
longer conjectural.
  \end{say}

  \begin{exmp}[Moduli of stable surfaces of general type]
\label{general type}

It is well known that the G.I.T. construction of a complete and projective
moduli space $\overline {\mathcal M}_{g}$ of stable curves does not work for
surfaces of general type.
As a bare minimum, the complete moduli space should
parameterize surfaces appearing as semistable degenerations of
nonsingular ones, and these have semi-log canonical singularities. In
particular, they can have singularities of arbitrarily large
multiplicities. On the other hand, by \cite{Mumford77} 3.19 a normal
surface with a singularity of multiplicity $\ge7$ is not
asymptotically Chow-(Hilbert-)stable,
so the usual G.I.T. construction does not go through.
 The following is a rather
general theorem which allows to circumvent this difficulty.
  \end{exmp}

    \begin{thm}[Koll\'ar, \cite{Kollar90}, 2.6]
Let $\mathcal C$ be an open class of $\mathbb Q$-polarized varieties with
Hilbert function $H(t)$. Assume that the corresponding moduli functor
${\mathcal M \mathcal C}$
is separated, functorially polarized, semipositive, bounded, complete
and has tame automorphisms. Then $\mathcal M \mathcal C$ is coarsely
represented by a projective scheme $\bold M \bold C$.
    \end{thm}

    \begin{say}
Applying this theorem to surfaces of general type one considers
stable surfaces, introduced for this purpose by Koll\'ar and
Shepperd-Barron in \cite{KollarShepherdBarron88}.
These are defined as surfaces with semi-log
canonical singularities (in particular, they are reduced but not
necessarily irreducible) such that $(\omega_{X}^{N})^{**}$ is ample
for some $N>0$.
Some of the conditions in the theorem are
easy to check. Separateness follows directly from the uniqueness of the
canonical model (in any dimension). Stable surfaces have finite
automorphism groups by Iitaka \cite{Iitaka82} 11.12
 (and this is also true in higher
dimensions). Completeness follows from the semi-stable reduction
theorem and the log Minimal Model Theory (in dimension $\dim X+1$, or,
as some say, $\dim X+1/2$). Semipositiveness is harder but it is proved in
\cite{Kollar90}. The boundedness is exactly what we are concerned with in
this paper. So if we want to include this application, the boundedness
theorem for surfaces with positive $K$ and constant $K^{2}$ should be
powerful enough as to include the case of semi-log canonical
singularities. And it is indeed, see~\ref{K2=const for semi-log
canonical} with empty $\mathcal C$.
    \end{say}

    \begin{say}
In fact,~\ref{K2=const for semi-log canonical} is strong enough so that
we are able to prove the projectiveness of the moduli space
$\overline {\mathcal M}_{(K+B)^{2}}$
of
surfaces with semi-log canonical and ample $K+B$ where $B$ is a reduced
divisor (or, more generally, an $\mathbb R$-divisor with coefficients in
a D.C.C. set $\mathcal C$) and $(K+B)^{2}=C$. This is a direct
generalization of the moduli space $\overline {\mathcal M}_{g,k}$ of pointed
stable curves (see \cite{Knudsen83}). This includes, for example, a
projective moduli space, an open subset of which parameterizes smooth K3
surfaces $X$ with ample and reduced divisors $B$ with normal crossings.
It is quite
interesting to see what is the relation between this moduli space and the
usual moduli space of polarized K3 surfaces. This
will be carried out in detail elsewhere.
    \end{say}

    \begin{say}
The following easy trick allows to reduce the boundedness of semi-log
canonical surfaces that are {\em a priori \/} not irreducible to the
boundedness of irreducible normal log canonical surfaces (see \cite{Kollar92}).
Let $X^{\nu}=\cup X_{m}$
be a decomposition into irreducible components as in
{}~\ref{reduce semi-log canonical to log canonical}.
Then
      \begin{displaymath}
K_{X}^{2}=  \sum (K_{X_{m}}+ B_{m})^{2}
      \end{displaymath}

So if we have enough information about possible values of
$(K+B)^{2}$ and know the
boundedness for normal surfaces with log canonical singularities, this
should help us in the general situation.
This  shows that it is natural to consider not
only the canonical divisor $K$ but also the canonical divisor with a
``boundary'' $K+B$.
The following example explains the importance of coefficients $b_{j}$ in
$B=\sum b_{j}B_{j}$ other than 1.
    \end{say}

    \begin{exmp}[Automorphisms of log surfaces of general type]
\label{cor auts of surfaces of general type}
The characteristic of the base field is assumed to be zero for this
application.

Consider a nonsingular surface with ample $K_{X}$. A general fact is
that the group of biregular automorphisms
of $X$ $\quad G=\Aut (X)$ is finite. Let $\pi:X\to Y$ denote a quotient
morphism. Then by the Hurwitz formula one has

      \begin{displaymath}
K_{X}= \pi^{*} (K_{Y}+ \sum (1-1/n_{j})B_{j})
      \end{displaymath}
where $B_{j}$ are the branching divisors. It follows that
$K_{Y}+B=K_{Y}+ \sum (1-1/n_{j})B_{j}$ is log terminal and that

      \begin{displaymath}
|\Aut(X)|= \frac{K_{X}^{2}}{(K_{Y}+B)^{2}},
      \end{displaymath}
so if $(K_{Y}+B)^{2})^{2} \ge 1/c$ then $|\Aut(X)|\le cK_{X}^{2}$. The
reader will certainly recognize that for curves this is the original
construction of Hurwitz which gives $\Aut(X)\le 42(2g_{X}-2)$ because
$2g_{Y}-2+\sum (1-1/n_{j})\ge 1/42$ if it is $>0$. It was Xiao's idea
to use the same construction for surfaces of general type in
\cite{Xiao91}. Using this and other methods he proves that for
nonsingular surfaces with ample canonical class one has
$|\Aut(X)|\le42^{2}K_{X}^{2}$. As an application of~\ref{DCC for
semi-log canonical} we have the following theorem.
    \end{exmp}

    \begin{thm}
\label{my thm on auts}
 Fix a D.C.C. set $\mathcal C$. Then for every surface $X$ with $K+B$ ample,
 semi-log canonical, and with $b_{j}\in\mathcal C$ one has the following
 bound for the
 automorphism group interchanging components of $B$ with the same
 coefficients
    \begin{displaymath}
|\Aut(X,B)|\le c(\mathcal C)(K+B)^{2}
    \end{displaymath}
 where the constant $c(\mathcal C)$ depends only on the set $\mathcal C$.
    \end{thm}
  \begin{proof}
The group
$G=\Aut(X,B)$ is known to be finite, cf.~\cite{Iitaka82}, 11.12.
Consider $\pi:X\to Y=X/G$.
Now use the same formulas
    \begin{eqnarray*}
K+B=\pi^{*}(K_{Y}+D) \text{ and }
|\Aut(X,B)|= \frac{(K+B)^{2}}{(K_{Y}+D)^{2}}
    \end{eqnarray*}
It easily follows that $K_{Y}+D$ is also semi-log canonical. The
coefficients of $D$ belong
to the set
    \begin{displaymath}
\{ 0\le 1-(1-\sum n_{j }b_{j})/m_{i}\le1  \,|\,  b_{j}\in\mathcal C,
\text{ } m_{i},n_{j}\in \mathbb N    \}
    \end{displaymath}
which is obviously also a D.C.C. set. Now we only need to see that a
D.C.C. set of positive numbers is bounded from below by a positive
constant.
  \end{proof}

  \begin{say}
Here is one more application of~\ref{DCC for semi-log canonical}
due to J.Koll\'ar  \cite{Kollar92}.
  \end{say}

  \begin{cor}[Uniform plurigenera of elliptic 3-folds]
For a smooth 3-fold of Kodaira dimension 2 in characteristic 0, there
exists an absolute constant $N$ such that $h^{0}(NK)\ne0$.
  \end{cor}

  \begin{say}
This was previously known for Kodaira dimension 0 (Kawamata,
\cite{Kawamata86}),
dimension 1 (Mori, \cite{Mori92})
and for dimension 3 and $\chi(\mathcal O_{X})\le M$ (Fletcher, \cite{Fletcher87}).

The  D.C.C. set used in this application is
  \begin{displaymath}
\mathcal C=\{ \frac{1}{12},\ldots{},\frac{11}{12},1-\frac{1}{k} \,|\,
k\in\mathbb N \}
  \end{displaymath}
  \end{say}

  \begin{say}
The following series of examples show that the restrictions of
{}~\ref{proved here} are
in a sense the weakest possible, and that none of them can be weakened
further.
  \end{say}

  \begin{exmp}
Consider the cone over an elliptic curve embedded
by a complete linear system of arbitrary degree in some projective
space.
The family of these cones is, evidently, not bounded. The
singularities are simple elliptic and they are log canonical, $-K$ is ample.
  \end{exmp}

  \begin{exmp}
There are infinitely many, I would say, hopelessly many
types of log Del Pezzo surfaces, i.e.\
surfaces with ample $-K$ and log terminal (=quotient in dimension 2
and characteristic 0) singularities. These include all surfaces $\mathbb
P^{2}/G$, $G$ a finite subgroup in $PGL(3)$, for example. So, for the
ample $-K$ the condition on singularities should be stronger than just
log canonical or log terminal. The $\varepsilon$-log terminal
condition seems to be the best substitute.
  \end{exmp}

  \begin{exmp}
The next example shows that
the condition $b_{j}\le 1-\varepsilon$ is necessary even for
smooth surfaces.
 Consider $\mathbb P^{2}$ with two lines
$B_{1}$ and $B_{2}$ and the surface obtained by
blowing up the point of the
intersection of these lines, then several times the point of the
intersection of the first line and the exceptional divisor.  Choosing
$b_{1}$ very close to 1 and taking for $K+B$ the full preimage
$f^{*}(K_{P^{2}}+b_{1}B_{1}+b_{2}B_{2})$
changed a little bit in the exceptional
divisors, one easily obtains an infinite sequence
of smooth surfaces with ample $-(K+B)$ and a strictly increasing Picard number.
The same construction works for surfaces with ample $K+B$ as well.
  \end{exmp}

  \begin{exmp}
This example shows that the set $\mathcal C$ in
{}~\ref{proved in this paper} has to satisfy the D.C.C.
 To
see this, take a rational ruled surface $\mathbb F_{e}$, $e\ge 2$ and
$B=(2+a)/4 (B_{1}+B_{2}+B_{3}+B_{4})$, where $B_{i}$ are general
elements in the linear system $|s_{e}+ef|$, $s_{e}$ is the exceptional
section, $f$ is a fiber. If $a$ is the positive root of the quadratic
equation $ea^{2}+(e-2)a=1$ then $K+B$ is ample and $(K+B)^{2}=1$. Note
that as $e\to \infty$, $(2+a)/4$ approaches its limit $ 1/2$ from
above.
  \end{exmp}

  \begin{rem}
It was conjectured by  V.V. Shokurov in \cite{Shokurov88} that all
sets naturally appearing in conjunction with log Minimal Model Program
have to satisfy either
 the descending or the ascending chain conditions. The first
example of this general phenomenon was given in  \cite{Alexeev89} for
the set of Fano indices of log Del Pezzo surfaces. See
\cite{Alexeev91a}, \cite{Alexeev93a},  \cite{GrassiKollar92}
for further examples.
  \end{rem}

\section{Some methods for proving boundedness}
\label{Some methods for proving boundedness}

  \begin{defn}
One says that a certain class of schemes $\mathcal B$ is bounded
if there exists a morphism $f:\mathcal X \to \mathcal S$ between two schemes
of finite type such that every scheme in $\mathcal B$ appears as one of
the geometric fibers of $f$, not necessarily in a one-to-one way. We do not
require that every geometric fiber of $f$ belongs to $\mathcal B$.
Usually, though, the class is defined by a combination of
algebraic conditions
some of which are open and others are closed. In this case one can
find a constructible subset in $S$, points of which parameterize
exactly elements of the class $\mathcal B$.
  \end{defn}

  \begin{defn}
In the same way, we say that a class $\mathcal B$ of schemes with closed
subschemes $\{(X,Z)\}$ is bounded if there exist three schemes of
finite type -- $\mathcal X$, a closed subscheme $\mathcal Z \subset \mathcal X$,
and $\mathcal S$, and a morphism $f:\mathcal X\to \mathcal S$ such that every
element of $B$ appears as a fiber of $f$.
  \end{defn}

  \begin{defn}
Finally, we can assign certain coefficients to subschemes of $X$
and then by boundedness of $\{(X,\sum b_{j}B_{j})\}$ we mean that
all $\{(X,B_{j})\}$ are bounded in the previous sense and, in
addition, that there are
only finitely many
possibilities for the sets of coefficients $\{b_{j}\}$.
  \end{defn}

  \begin{defn}
A polarization on a scheme $H$ is a class of an ample Cartier divisor.
A $\mathbb
Q$-polarization on a normal variety is a $\mathbb Q-$Cartier divisor a,
positive multiple of which is a polarization.
It is possible to define an $\mathbb R$-polarization on some nonnormal
varieties too
when there is a suitable notion of an $\mathbb R$-divisor. A semi-log
canonical scheme with an ample $\mathbb R$-Cartier divisor $K+B$ would be
an example.
  \end{defn}

  \begin{say}
Consider a class of polarized reduced schemes over an algebraically
closed  field $k$ $(X,H)$ with a fixed
Hilbert function $\mathcal H(t)=\chi(tH)$. It is known that this class is
bounded provided any of the following conditions are satisfied:

    \begin{enumerate}
      \item $\dim X=2$ (Matsusaka \cite{Matsusaka86} for normal
      surfaces, Koll\'ar \cite{Kollar85} for the general case).
      \item $\dim X=3$, $X$ are normal and $\Char k=0$ ( \cite{Kollar85}).
      \item $X$ are nonsingular and $\Char k=0$ (Matsusaka's Big
      Theorem, see \cite{Matsusaka86}).
    \end{enumerate}
  \end{say}

  \begin{say}
Moreover, if $\Char k=0$ and $X$ are normal, then only the first two
coefficients
of $\mathcal H(t)$, i.e.\ , up to constants, $H^{\dim X}$ and $H^{\dim
X-1}K_{X}$, are important, by the Riemann-Roch inequalities of
Koll\'ar-Matsusaka
\cite{KollarMatsusaka83}.
In dimension two this is also true in arbitrary characteristic, see
Lemma 2.5.2 \cite{Kollar85}.
  \end{say}

  \begin{lem}\label{main method for proving boundedness}

    \begin{enumerate}
      \item Fix $C>0$ and $N\in \mathbb N$, and
            let $\mathcal B= \{ (X,H) \}$ be a class of normal $\mathbb Q$-polarized
            surfaces such that $NH$ is Cartier, $H^{2}\le C$.
            Then
            the class $\mathcal B$ is bounded.
      \item Fix $C,C'>0$ and $N\in \mathbb N$. Let
            $\mathcal B= \{ (X,H,B) \}$ be a class of normal
            $\mathbb Q$-polarized surfaces with
            subschemes, where $B$ is either a divisor with $HB\le C'$ or
            a $0$-dimensional subscheme of length $\le$ $C'$. Then
            the class $B$ is bounded.
    \end{enumerate}
  \end{lem}
  \begin{proof}
The first part is just a reformulation of the above remarks. In $
(2)$, since the class $\{(X,H)\}$ is bounded, all the surfaces can be
embedded by a uniform multiple of $H$ in the same projective space
$\mathbb P$. Then the subschemes $B$ satisfying the conditions above are
parameterized by a certain Hilbert scheme $\bold H$, and the pairs
$(B\subset X)$ are parameterized by a closed subscheme of $\mathbb P\times \bold
H$.
  \end{proof}

  \begin{lem}
\label{boundedness by blowingup}
 Let $\{X\}$ be a certain class of schemes and assume that every scheme
$X$ is isomorphic to a blowup $Bl_{Z}Y$ of a scheme $Y$ at a subscheme
$Z$. If the class $\{(Y,Z)\}$ is bounded then the class $\{X\}$ is
bounded.
  \end{lem}
  \begin{proof}
After subdividing $\mathcal S$ we can assume that both $\mathcal Y$ and $\mathcal
Z$ are flat over $\mathcal S$. Then
 $Bl_{\mathcal Z}\mathcal Y$ will give a required family.
  \end{proof}

  \begin{say}
I shall use one more method for proving boundedness given below, and a
similar method for proving the descending chain condition.
They sound almost
trivial but I still would like to formulate them explicitly.
  \end{say}

  \begin{thm}
\label{kinky method for proving boundedness}
Let  $\mathcal B$ be a certain class of schemes. Assume that for every
infinite sequence $\{X_{s}\in \mathcal B\}$ there exists an infinite
subsequence $\{ X_{s_{k}} \}$ which is bounded. Then the class $\mathcal
B$ is bounded.
  \end{thm}

  \begin{thm}
\label{kinky method for proving DCC}
Let $\mathcal C$ be an ordered set. Assume that for every
infinite sequence $\{ x_{s}\in \mathcal C\}$ there exists an infinite
nondecreasing
subsequence $\{ x_{s_{k}} \}$. Then the set $\mathcal
C$ satisfies the descending chain condition.
  \end{thm}

\section{Additional definitions and easy technical results}
\label{Additional definitions and easy technical results}

  \begin{defn}
    For two $\mathbb R$-divisors we write $D_{1}\ge D_{2}$
    (resp.~$D_{1}> D_{2}$) if
    $D_{1}-D_{2}$ is effective, i.e.\  the coefficients for all prime
    components are nonnegative
 (resp. is effective and nonzero).
  \end{defn}

  \begin{defn}
    For two $\mathbb R$-divisors on possibly different normal varieties
of the same dimension we write
    $D_{1}\underset{c}{\ge} D_{2}$ (resp.
    $D_{1}\underset{c}{>} D_{2}$) if
    $$\underset{n\to\infty}{\varliminf}
    \frac{h^{0}(nD_{1})-h^{0}(nD_{2})}{n^{\dim}}\ge0$$
    (resp. is strictly positive). Here $n$ is assumed to be divisible
    enough.
  \end{defn}

  \begin{defn}
One says that an $\mathbb R$-Cartier divisor $D$ is {\em nef\/} if
$DC\ge0$ for any effective curve $C$.
  \end{defn}

  \begin{defn}
The {\em Iitaka dimension} or the {\em $D$-dimension} of a Weil divisor
$D$ on a normal variety,
denoted
by $\nu(D)$, is
the number  in the formula $h^{0}(nD) \sim  n^{\nu(D)}$ as $n\to\infty$,
unless $h^{0}(nD)=0$ for all $n>0$, in which case $\nu(D)$
is defined to be equal to $-\infty$. A divisor with $\nu(D)=\dim X$ is
called {\em big}. The {\em Kodaira dimension} of a variety $X$ is
defined as $\nu(K_{Y})$ where $Y$ is a resolution of singularities of $X$.
  \end{defn}

  \begin{lem}
\label{easy inequalities}
      \begin{enumerate}
        \item $D_{1}\ge D_{2}$ implies $D_{1}\underset{c}{\ge} D_{2}$,
        \item if $H^{i}(nD_{1})=H^{i}(nD_{2})=0$ for $i>0$, $n\gg0$,
              then $D_{1}\underset{c}{\ge} D_{2}$ (resp.
              $D_{1}\underset{c}{>} D_{2}$) is equivalent to
              $D_{1}^{\dim}\ge D_{2}^{\dim}$ (resp.
              $D_{1}^{\dim}>D_{2}^{\dim}$),
        \item in particular, (2) is applicable when $D_{1}$ and
        $D_{2}$ are ample.
      \end{enumerate}
  \end{lem}
  \begin{proof}
    (1) is clear, (2) and (3) trivially follow from the Riemann--Roch
    formula and  the Serre  vanishing theorem.
  \end{proof}

  \begin{lem}
\label{inequalities}
    Let $D_{1}$ be an ample $\mathbb R$-divisor on a normal variety $X$, $f:Y\to
    X$ be a birational projective morphism with $Y$  normal,
    and write
    $$D_{2}=f^{*}(D_{1})+\sum e_{j}E_{j}+\sum f_{i}F_{i}$$
    where the divisors $F_{i}$ are exceptional for $f$ and $E_{j}$ are
    not.
    Then
    \begin{enumerate}
      \item if $e_{j}\le0$ then $D_{2}\underset{c}{\le}D_{1}$
      \item if, moreover, for some index  $e_{j_{0}}<0$ or $f_{i_{0}}<0$ then
            $D_{2}\underset{c}{<}D_{1}$
      \item assume that
	    $e_{j}\le0$,
	    $D_{2}$ is {\em nef\/};
	    then all $f_j\le0$ i.e.\  $D_{2}\le f^{*}(D_{1})$
      \item assume that
	    $e_{j}\le0$,
	    $D_{2}$ is {\em nef\/}
            and that $D_{1}\underset{c}{=}f^{*}(D_{2})$; then
            $D_{2}=D_{1}$ i.e.\  all $e_{j}=f_{i}=0$
    \end{enumerate}
  \end{lem}

  \begin{proof}
   Most of the statements are elementary, others follow from the
   Negativity of Contractions Lemma, see \cite{Shokurov91}, 1.1 or
   \cite{FAAT}, 2.19.
  \end{proof}

  \begin{lem}
\label{only surfaces}
    Let $D_{1}$ be a divisor on a nonsingular surface $X$ such that
    $H^{i}(nD_{1})=0$ for $i>0,$ $n\gg0$. Let $D_{2}$ be a nonzero
    divisor such that $-D_{2}$ is not quasieffective
(i.e.\  there exists an ample divisor $H$ such
that  $-D_{2}H<0$). Assume that
    $D_{2}\underset{c}{\le}D_{1}$. Then $D_{2}^{2}\le D_{1}^{2}$.
  \end{lem}
  \begin{proof}
    Use the Riemann--Roch theorem for $nD_{1},nD_{2}$
and the fact that $h^{2}(nD_{2})=h^{0}(K-nD_{2})=0$
for $n\gg0$.
  \end{proof}

\section{The diagram method}
\label{The diagram method}

  \begin{say}
In proving the boundedness results for surfaces with nef $-(K+B)$ the
main part will be to bound rank of the Picard group. A way for
doing this, called the diagram method,
 comes from the theory of reflection groups in hyperbolic
spaces. Given a convex polyhedron in a hyperbolic space it is possible
to bound the dimension of this space purely combinatorially, provided
subpolyhedra of the main polyhedron
satisfy certain arithmetic properties. The diagram method was
successfully applied to surfaces in the works of V.V.Nikulin
 \cite{Nikulin90c}, \cite{Nikulin89b}, \cite{Nikulin90a},
\cite{Nikulin89a}
and of the author \cite{Alexeev89}, and also to Fano 3-folds in
 \cite{Nikulin90b}, \cite{Nikulin93a}.
The vector space in question is $\Pic(X)\otimes \mathbb R$ which is
hyperbolic by the Hodge Index theorem or some linear subspace of it,
the polyhedron is generated by exceptional
curves.

In the next section I shall give one more application of the diagram method.
But first we need a few definitions and facts.
  \end{say}

  \begin{defn}
An {\em exceptional curve\/} on a surface is an irreducible curve $F$
with $F^{2}<0$. The set of all exceptional curves on the surface $X$
will be denoted by $\Exc(X)$. A {\em $(-n)$-curve} is a smooth rational
curve $F$ with $F^{2}=-n$.
  \end{defn}

  \begin{say}
To a set of exceptional curves we  associate a {\em
weighted graph}. Each curve $F$ corresponds to a vertex of weight
$-F^{2}$ and two vertices $F_{1}$ and $F_{2}$ are connected by an
edge of weight  $F_{1}F_{2}$. An edge of weight 1 will be called {\em
simple}. We also  assign to every vertex a nonnegative
number,
the arithmetical genus of the corresponding curve. However, in the
situation of interest to us all exceptional curves will have genus zero, as
{}~\ref{genera only zero} shows.
The (possibly infinite) graph corresponding to all exceptional curves
on a surface $X$ will be denoted by $\Gamma(\Exc(X))$.
  \end{say}

  \begin{defn}
A finite set of exceptional curves $\{ F_{i} \,|\, i=1\ldots{}r  \}$
(and the corresponding weighted
graph) is called {\em elliptic\/} (resp.~ {\em parabolic, hyperbolic\/}) if the
matrix $(F_{i_{1}}F_{i_{2}})$ has signature $(0,r)$
(resp.~$(0,r-1)$, $(1,r-1))$.
  \end{defn}

  \begin{defn}
A finite set of exceptional curves
(and the corresponding weighted
graph) is called {\em Lanner\/} if it is hyperbolic but any proper subset
of it is not.
  \end{defn}

  \begin{say}
The following theorem will be the most important for our purposes, see
 \cite{Nikulin89a}, 3.4--6 for the proof.
  \end{say}

  \begin{thm}
\label{Nikulin's main theorem}
Let $X$ be a nonsingular surface such that the Iitaka dimension of the
anticanonical divisor is nonnegative. Let us assume that for certain
constants $d$,$c_{1}$,$c_{2}$ the following conditions hold:
    \begin{enumerate}
      \item the diameter of any Lanner subgraph $\mathcal L\subset\Gamma(\Exc(X))$
            does not exceed $d$;
      \item if $\nu(-K)=2$ then for any connected elliptic subgraph
            $\mathcal E \subset \Gamma(\Exc(X))$ with $n$ vertices, the number
            of (unordered) pairs of its vertices on distance $\rho$,
            where $1\le \rho \le d-1$, does not exceed $c_{1}n$, and the
            number of pairs on distance $\rho$, where $d\le \rho \le
            2d-1$, does not exceed $c_{2}n$;
      \item if $\nu(-K)=1$ then for any $(-1)$-curve $E$ for which
            $EP>0$ ($P$ is the positive part of the Zariski
            decomposition for $-K$ and won't be used later) and
            for any connected elliptic subgraph
            $\mathcal E \subset \Gamma(\Exc(X))$ with $n+1$ vertices
            containing $E$, the number
            of pairs
            of its vertices
            different from $E$ and on distance $\rho$,
            where $1\le \rho \le d-1$, does not exceed $c_{1}n$, and the
            number of pairs on distance $\rho$, where $d\le \rho \le
            2d-1$, does not exceed $c_{2}n$;
      \item if $\nu(-K)=0$, $-K=\sum_{b_{j}>0}b_{j}B_{j}$
            then for any connected elliptic subgraph
            $\mathcal E \subset \Gamma(\Exc(X))$ with $n+m$ vertices,
            $m$ of which $E_{1}\ldots{}E_{m}$
            correspond to $(-1)$- or $(-2)$-curves different from
            $B_{j}$, the number
            of pairs of its vertices different from
            $E_{1}\ldots{}E_{m}$ and on distance $\rho$,
            where $1\le \rho \le d-1$, does not exceed $c_{1}n$, and the
            number of pairs on distance $\rho$, where $d\le \rho \le
            2d-1$, does not exceed $c_{2}n$.
    \end{enumerate}
Then
    \begin{enumerate}
      \item if $\nu(-K)=2$ then $\rk\Pic(X)\le96(c_{1}+c_{2}/3)+69$;
      \item if $\nu(-K)=1$ then $\rk\Pic(X)\le96(c_{1}+c_{2}/3)+70$;
      \item if $\nu(-K)=0$ then $\#(B_{j})\le96(c_{1}+c_{2}/3)+68$.
    \end{enumerate}
  \end{thm}

  \begin{say}
We  need a few definitions for blowing up and down weighted graphs. They
are merely reformulations on the language of graphs of
usual operations of blowing up points on a nonsingular surface.
  \end{say}

  \begin{defn}
A weighted graph is said to be {\em minimal} if it does
not contain vertices of weight 1 and arithmetical genus 0.
  \end{defn}

  \begin{defn}
{\em Blowing up} a vertex $F$ is the operation on  a weighted
graph consisting
of adding a new vertex $E$ of weight 1
and of arithmetical genus 0,
connected only with the vertex
$F$ by a simple edge and increasing weight of $F$ by 1.
  \end{defn}

  \begin{defn}
{\em Blowing up} a simple edge $F_{1}F_{2}$ is the operation on a
weighted graph consisting
of adding a new vertex $E$ of weight 1
and of arithmetical genus 0,
 connected only with the vertices
$F_{1}$ and $F_{2}$ by simple edges, removing the edge between them
and increasing weights of $F_{1}$ and $F_{2}$ by 1.
  \end{defn}

  \begin{say}
 Blowing up
can be easily defined in a more general situation but we won't need
it.
  \end{say}

  \begin{defn}
{\em Blowing down} is the inverse operation  to blowing up.
  \end{defn}

  \begin{defn} {\em The canonical class} $K=K(\Gamma)$ of a weighted
    graph $\Gamma$ is the function on vertices defined by the formula
    \begin{displaymath}
KF_{i}=-F_{i}^{2}-2+2p_{a}(F_{i})
    \end{displaymath}
  \end{defn}

  \begin{defn} {\em The log discrepancies} $f_{i}$ for a finite
    weighted graph $\Gamma$ are defined as solutions of the following
    system of linear equations (if exist):
    \begin{displaymath}
(K+\sum (1-f_{i})F_{i})F_{j}=0 \text{\quad for all }j
    \end{displaymath}
Note that this system has a unique solution if the matrix
$(F_{i}F_{j})$ is invertible (for example if the graph $\Gamma$
is elliptic or hyperbolic).
  \end{defn}

  \begin{defn}
A weighted graph $\Gamma$ is said to be {\em log terminal} if for every
elliptic subgraph $\Gamma' \subset \Gamma$ all log discrepancies of
$\Gamma'$ are positive.
  \end{defn}

  \begin{defn}
\label{def of condition star for graphs}
We say that a finite graph $\Gamma=\{F_{i}\}$ satisfies {\em the
condition $*(\varepsilon)$} if there exist constants $0\le b_{i}\le
1-\varepsilon<1$ such that
$(K+\sum b_{i}F_{i})F_{j}\le0$ for all vertices $F_{j}$.
  \end{defn}

  \begin{lem}
\label{condition star}
If $\Gamma$ satisfies $*(\varepsilon)$ then every subgraph
$\Gamma_{1} \subset \Gamma$ and every graph $\Gamma_{2}$ obtained from
$\Gamma$ by blowing down several vertices of weight 1 also satisfy
$*(\varepsilon)$.
  \end{lem}
  \begin{proof}
Evident.
  \end{proof}

  \begin{lem}
If an elliptic graph
$\Gamma$ satisfies $*(\varepsilon)$ then for all the log
discrepancies $f_{i}\ge 1-b_{i} \ge \varepsilon$.
  \end{lem}
  \begin{proof}
This is well known and follows easily from the negative definiteness
of $(F_{i}F_{j})$, see for example \cite{Alexeev92}, 3.1.3.
  \end{proof}

  \begin{lem}
\label{genera only zero}
If an elliptic graph
$\Gamma$ satisfies $*(\varepsilon)$ then
every vertex $F$ has arithmetical genus 0 and
its
weight does not exceed $2/\varepsilon$.
  \end{lem}
  \begin{proof} Follows from

    \begin{eqnarray*}
-2 \le  2p_{a}(F)-2 = (K+F)F = (K+ (1-\varepsilon)F)F + \varepsilon
F^{2} \le \\
(K+\sum b_{j}B_{j})F + \varepsilon F^{2} \le \varepsilon F^{2} <0
    \end{eqnarray*}
\renewcommand{\qed}{}
\qed
  \end{proof}

  \begin{thm}
\label{minimal elliptic graph}
Every minimal elliptic log terminal
graph is a
tree with at most one fork and of type $A_{n}$, $D_{n}$ or $E_{6,7
\text{ or }8}$. There exists a constant $S_{1}(\varepsilon)$ depending
only on $\varepsilon$ such that $\sum (-F^{2}_{i}-2) \le
S_{1}(\varepsilon)$
if such  a graph satisfies $*(\varepsilon)$.
  \end{thm}
  \begin{proof}
The first part of the statement is well known,
see for example \cite{Alexeev92}. The second part follows
from the explicit description of elliptic graphs with all log
discrepancies $\ge\varepsilon>0$ given in \cite{Alexeev93a}, 3.3.
  \end{proof}

  \begin{thm}
\label{minimal Lanner graph}
There exist $\le 14(2/\varepsilon)+29$
Lanner graphs with simple edges
such that every Lanner graph
with more than 5 vertices satisfying $*(\varepsilon)$ can be
obtained from one of them by blowing up several vertices and edges.
Each of these graphs is a tree or a cycle or a cycle and one more
vertex.
Each of these graphs has only simple edges and
every vertex has at most 3 neighbors.
  \end{thm}
  \begin{proof}
It follows from the theorems of V.V.Nikulin, \cite{Nikulin89a}, 4.4.18,19,21
which are valid for arbitrary log terminal Lanner graphs.
Some of the graphs in \cite{Nikulin89a} have a vertex of arbitrary
positive weight $b$ but in our situation $b\le 2/\varepsilon$ by
{}~\ref{condition star} and~\ref{genera only zero}.
  \end{proof}

  \begin{say}
A typical example of a Lanner graph in the above statement is the
chain containing three vertices of weights 1, 1 and $b\ge 1$. We won't
need explicit description of these graphs, just knowing that for every
$\varepsilon$ there are finitely many of them will be sufficient.
  \end{say}

\section[Boundedness for surfaces with nef $-(K+B)$
]{Boundedness for surfaces with nef $\bold -(K+B)$}
\label{Boundedness for surfaces with nef -(K+B)}

  \begin{lem}
\label{graphs on good surfaces satisfy the condition star}
Let $X$ be a nonsingular surface and assume that $-(K+\sum
b_{j}B_{j})$ is nef, where $0\le b_{j} \le 1-\varepsilon <1$.
Then
every finite subgraph
$\Gamma\subset \Gamma(\Exc(X))$ satisfies the condition $*(\varepsilon)$.
  \end{lem}
  \begin{proof}
Evident.
  \end{proof}

  \begin{lem}
\label{structure when -K ge0}
Let $X$ be a nonsingular surface and assume that $-(K+\sum
b_{j}B_{j})$ is nef, where $0\le b_{j} \le 1-\varepsilon <1$.
Then one of the following is true:
    \begin{enumerate}
      \item all $b_{j}=0$ and $K$ is numerically trivial
      \item $X$ is a rational surface, obtained by blowing up several points
            from $\mathbb P^{2}$ or $\mathbb F_{n}$ with $n\le 2/\varepsilon$
      \item $X$ is an elliptic ruled surface without exceptional
            curves, i.e.\
            a projectivization of a rank 2 locally free
            sheaf $\mathcal E$ on an elliptic curve $C$ and $\mathcal E$ is isomorphic
            to $\mathcal O\oplus \mathcal F$, $\mathcal F\in \Pic^{0}(C)$
            or to one of the only two, up to
            tensoring with an invertible sheaf,
            nonsplilttable rank 2 bundles on $C$
    \end{enumerate}
  \end{lem}
  \begin{proof}
This lemma is practically proved in \cite{Nikulin89a}, 4.2.1
under weaker assumptions on $b_{j}$ and in arbitrary
characteristic. We only need to see that in the case (2) one has $n\le
2/\varepsilon$ by~\ref{genera only zero}.
  \end{proof}

  \begin{thm}
\label{thm A1}
Let $X$ be a nonsingular surface and assume that $-(K+\sum
b_{j}B_{j})$ is nef, and $0\le b_{j} \le 1-\varepsilon <1$. Then
there exists a constant $A_{1}(\varepsilon)$ which depends only on
$\varepsilon$ such that
    \begin{displaymath}
\rk\Pic(X)\le A_{1}(\varepsilon)
    \end{displaymath}
  \end{thm}

  \begin{proof}
In order to prove this theorem we have to check the conditions (1) and
(2--4) for graphs satisfying $*(\varepsilon)$. Among them, (1) is the
hardest.
Given an arbitrary graph satisfying $*(\varepsilon)$ and containing a
vertex of weight 1, we  can blow this vertex down and the new graph will
satisfy the same condition by~\ref{condition star}. After blowing down
several vertices of an elliptic graph we get a minimal graph
described in~\ref{minimal elliptic graph}.
It is well known that for an elliptic graph
blowing down vertices of weights 1 in any order yields the same graph,
so we can
backtrack the situation by blowing up edges and vertices on the minimal
elliptic graph {\em in arbitrary order}. All intermediate graphs
again should satisfy $*(\varepsilon)$.

For a Lanner graph, removing any two vertices gives an elliptic graph.
So, basically, we can do the same thing. Contracting several vertices
we obtain one of the finitely many graphs of
{}~\ref{minimal Lanner graph}, and we can backtrack
the situation blowing up,
in any  order,
edges and vertices that do not
affect two
arbitrarily chosen vertices.
Again, all intermediate graphs should satisfy
$*(\varepsilon)$, moreover, they all should be Lanner.

Using these considerations, the proof easily follows from the lemmas
{}~\ref{graphs: blow up vertices}, \ref{graphs: few neighbors},
\ref{graphs: blow up edges}, \ref{graphs: blow up till get E9} below.

    \begin{lem}
\label{graphs: blow up vertices}
Fix a weighted graph $\Gamma$ and pick one of its vertices $F$.
Blow
it up to get the vertex $E_{1}$, then blow up $E_{1}$ to get $E_{2}$,
and so on. Call the intermediate graphs
$\Gamma_{1},\Gamma_{2}\ldots{}$
Then for $k\gg0$ the graph $\Gamma_{k}$ is not Lanner.
    \end{lem}
    \begin{proof}
It is easy to see, using only elementary linear algebra, that there
exists a fractional linear function $f(k)$
(which is an appropriately normalized determinant)
so that the condition for
the graph $\Gamma-E_{k}$ not to be hyperbolic
 is equivalent to $f(k)\ge0$. But then
the condition for $\Gamma$ to be hyperbolic is equivalent to
$f(\infty)<0$ and hence $f(k)<0$ for $k\gg0$.
    \end{proof}

    \begin{lem}
\label{graphs: few neighbors}
In any  satisfying $*(\varepsilon)$ graph
which is elliptic or is Lanner with more than 5  vertices,
each
vertex has at most $2/\varepsilon-2$ neighbors.
    \end{lem}
    \begin{proof}
Indeed, in each of the initial graphs in~\ref{minimal elliptic graph}
and~\ref{minimal Lanner graph} every vertex
has at most 3 neighbors,
so to produce a vertex with $d$ neighbors one has to
blow up some vertex $\ge d-3$ times, which means that the weight of
this vertex will be $\ge d-2$. Now use~\ref{genera only zero} .
    \end{proof}

    \begin{lem}
\label{graphs: blow up edges}
Fix an elliptic weighted graph $\Gamma$ consisting of two vertices
connected by a
simple edge.
Then among all graphs $\Gamma'$
obtained from $\Gamma$ by blowing up only edges,
there exist $\le S_{2}(\varepsilon)$ graphs satisfying $*(\varepsilon)$,
where the last function $S_{2}(\varepsilon)$ depends only on $\varepsilon$.
    \end{lem}
    \begin{proof}
For any vertex $F$ in any such graph $\Gamma'$ let us introduce {\em
the height} $h(F)$ as the minimal number of blowups needed to obtain
this vertex from $\Gamma$. Consider this minimal sequence of blowups.
Note that it is unique.
On each  step we get two chains of vertices of weight $\ge2$ and a
single vertex of weight 1 between them. After every blowup the sum
$\sum (-F^{2}_{i}-2)$ increases by 1. Therefore by~\ref{minimal
elliptic graph}
$h(F)\le S_{1}(\varepsilon)/2$ and there are $\le 2^{S_{1}/2-1}$ such
vertices. Finally,
there exist only finitely many graphs
$\Gamma'$ that contain only vertices of bounded height.
    \end{proof}

    \begin{lem}
\label{graphs: blow up till get E9}
Fix an elliptic weighted graph $\Gamma$ consisting of two vertices
connected by a
simple edge. Blow up the edge to get
the vertex $E_{1}$, then blow up $E_{1}$ to get $E_{2}$,
and so on. Call the intermediate graphs
$\Gamma_{1},\Gamma_{2}\ldots{}$.
Then for $k>5$ the graph $\Gamma_{k}$ is not log terminal, i.e.\ it
does not satisfy
$*(\varepsilon)$ for any $\varepsilon>0$.
    \end{lem}
    \begin{proof}
Indeed, $\Gamma_{6}-E_{6}$ has type $E_{9}$ or worse.
    \end{proof}

{\em End of the proof of~\ref{thm A1}.} The conditions (2--4) of
{}~\ref{Nikulin's main theorem} follow
easily from~\ref{graphs: few neighbors}.
Indeed, the number of pairs in (2--4) on distance
$\rho$ is bounded by $n/2(2/\varepsilon-2)^{\rho}$.

To prove that the condition (1) of
{}~\ref{Nikulin's main theorem} is satisfied, consider any of
the finitely many Lanner graphs of
{}~\ref{minimal Lanner graph}. We can assume that this graph has at
least 6 vertices. We shall prove that
there is a bound, in terms of $\varepsilon$, on the number of possible
blowups that can be done preserving the condition $*(\varepsilon)$ and
the property of the whole graph to remain Lanner.
As was mentioned before, fixing any two vertices, the order of
blowups not affecting these vertices is unimportant.

First, by
{}~\ref{graphs: blow up vertices},
vertices can be blown up only finitely many times
(again, there is a bound in terms of $\varepsilon$). Then an edge
between some two vertices
should be blown up. By~\ref{graphs: blow up edges},
different edges between these two vertices
can be blown up only finitely many times. Then some of the vertices
may acquire more neighbors, but again their number is limited by
{}~\ref{graphs: few neighbors}.
After that one of the new branches may grow longer, but its length is
bounded by~\ref{graphs: blow up till get E9}.
Then, again there may be a few edges blown up,
and then a few vertices may acquire some new neighbors. Finally, on
this stage if some of the branches grow longer, a minimal elliptic
intermediate
subgraph should appear which has more than one fork. But this is
impossible by~\ref{minimal elliptic graph}.

Therefore, all conditions (1) and (2--4) of
{}~\ref{Nikulin's main theorem} are satisfied with
constants
$d(\varepsilon)$, $c_{1}(\varepsilon)$, $c_{2}(\varepsilon)$ that depend only
on $\varepsilon$. In the cases $\nu(-K)=2\text{ or }1$ we immediately get
the upper bound on the Picard number. In the case $\nu(-K)=0$,
$-K=\sum_{b_{j}>0} b_{j}B_{j}$ one has
    \begin{displaymath}
K^{2}\ge -(\#B_{j})(2/\varepsilon)
    \end{displaymath}
 by~\ref{genera only zero}
and if $X$ is rational
we have $\rk\Pic(X)=10-K^{2}$ by Noether's formula. In the remaining two
cases of
{}~\ref{structure when -K ge0} one certainly has $\rk\Pic(X)\le20$.
  \end{proof}

  \begin{thm}
\label{bound for nonsingular when -K nef}
Fix $\varepsilon>0$. Consider all nonsingular surfaces
$X$ with an $\mathbb R$-divisor $B=\sum b_{j}B_{j}$ such that
$0\le b_{j}\le 1-\varepsilon <1$ and $-(K+B)$ nef excluding only those
for which at the same time $K_{X}$ is numerically trivial and $B$ is
zero.
Then the class $\{X\}$ is bounded.
  \end{thm}
  \begin{proof}
Indeed, in the case (3) of
{}~\ref{structure when -K ge0}
 there are only three deformation types.
All other surfaces, except $\mathbb P^{2}$,
are obtained from $\mathbb F_{n}$, $n\le 2/\varepsilon$
by $\le A_{1}(\varepsilon)-2$ blowups.
Now use~\ref{boundedness by blowingup}.
  \end{proof}

  \begin{thm}
\label{bound for general when -K nef}
  Fix $\varepsilon>0$.
  Consider all projective surfaces $X$
  with an $\mathbb R$-divisor $B=\sum b_{j}B_{j}$ such that
  $K_{X}+B$ is MR $\varepsilon$-log terminal and $-(K_{X}+B)$ is
  nef excluding only those for which at the same time $K_{X}$
  is numerically trivial, $B$ is zero and $X$ has at worst Du Val
  singularities.
  Then the class $\{X\}$ is bounded.
  \end{thm}
  \begin{proof}
By~\ref{bound for nonsingular when -K nef} we already know that the
class of minimal desingularizations $\{Y\}$ of surfaces $\{X\}$ is
bounded. This, however, does not yet guarantee the boundedness of the
class $\{X\}$. To prove it we shall use a ``sandwich'' argument: we
shall prove that there exist two birational morphisms $Y\to X\to Z$
with $\{Z\}$ also bounded.

By~\ref{structure when -K ge0}
 we can assume that $X$ is rational, moreover, by
{}~\ref{thm A1} the rank
of the Picard group of the minimal resolution of $X$ is effectively bounded.

We consider the following cases:
    \begin{case}$B$ is nonempty or $K$ is not numerically trivial.

We want to show that we can find a contraction $X\to Z$ so that
$-K_{Z}$ is ample or is relatively ample with respect to some $\mathbb
P^{1}$-fibration.
First, assume that there exists an $\mathbb R$-divisor $D$ such that
$K+B+D$ is
numerically trivial
and satisfies the same condition
$*(\varepsilon)$. According to our conditions, $B+D$ is nonempty.
If we decrease one of the coefficients in $K+B+D$, the result will be
not nef. Therefore, according to the Minimal Model Program there exists an
extremal ray and the corresponding contraction.
If $E$ is a component of $B+D$ with $E^{2}\le 0$
then $K+B+D-eE$ for $0<e\ll1$ has a nonnegative
intersection with $E$. Hence, there always exists an extremal contraction
which does not contract $E$,
so the image of $B+D$ will be again nonempty.

 Performing several such contractions we
arrive either at a surface $Z$ with $\rk\Pic(Z)=1$ and ample $-K_{Z}$
or at a surface $Z$ with $\rk\Pic(Z)=2$ which has a $\mathbb
P^{1}$-fibration. Because the rank of the Picard group of the minimal
resolution is bounded and the self-intersection numbers of exceptional
curves are bounded from below,
there are only finitely many types of singularities that $Z$ can have and
all of them can be effectively described in terms of $\varepsilon$.
Therefore, in the case of $\rk\Pic(Z)=1$
 a fixed multiple of $-K_{Z}$ is an ample Cartier divisor and
$K_{Z}^{2}$ is effectively bounded from above.
In the case of a $\mathbb P^{1}$-fibration,
$-K_{Z}+(\lceil 2/\varepsilon \rceil -1)\times fiber$ is ample, and
the same argument applies. The class $\{Z\}$ is bounded by
{}~\ref{main method for proving boundedness}.

Now, instead of proving that such a divisor $D$ as above exists, consider an
arbitrary ample divisor $A$ on $X$ and     observe that  we  can as
well use $D'/N$ in the place of $D$, where $D'$ is a general element
of a very ample linear system $|N(\delta A-K-B)|$ for $N\gg0$,
$0<\delta\ll  \varepsilon$.

The surface $Y$ is obtained from $Z$ by blowing up a certain closed
subscheme.
Again, from the boundedness of the Picard number of the resolution of
$Y$ and the squares of exceptional curves, we see that the length of
this subscheme is effectively bounded in terms of $\varepsilon$.
Now use~\ref{boundedness by blowingup}.
    \end{case}

    \begin{case}$B$ is empty, $K$ is numerically trivial but $X$ has
      worse than Du Val singularities.

Consider
a partial resolution $f:Y\to X$ dominated by
the minimal resolution of singularities of $X$ such that,
letting $K_{Y}+B^{Y}=f^{*}(K_{X})$, $B^{Y}$ has positive coefficients
in all exceptional divisors of $f$.
By the previous case,
there exists a polarization $H$ on $Y$ with $H^{2}$ and $K_{Y}H$ bounded.
Also, there exists
a lower bound, in terms of $\varepsilon$,
for $b_{j}^{Y}$. Since for any ample divisor $H$,
      \begin{displaymath}
\sum b_{j}^{Y}B_{j}H = -K_{Y}H
      \end{displaymath}
we can bound all $B_{j}H$ and, therefore,
for $H'=f_{*}(H)$
we can bound $H^{\prime 2}$, $H'K_{X}$ and an integer $N$ such that $NH$ is
Cartier. Then apply~\ref{main method for proving boundedness} again.
     \end{case}
\renewcommand{\qed}{}
  \end{proof}

  \begin{cor}
\label{bound for log terminal when -K nef}
The class of projective complex surfaces
$X$ with $K$ numerically trivial and  log terminal and
singularities worse than only Du Val, is bounded.
  \end{cor}
  \begin{proof}
Indeed, by \cite{Blache92}, Theorem C,
$GK_{X}$ is Cartier for some $G\in\{ 1,2\ldots{}21 \}$,
so every such surface is $1/21$-log terminal.
  \end{proof}

  \begin{thm}
\label{thm A2 for general}
Fix $\varepsilon>0$. Then there exists a constant $A_{2}(\varepsilon)$
such that for any projective surface
$X$ with $-(K+B)$  nef and $K+B$
MR $\varepsilon$-log terminal,
    \begin{displaymath}
\sum b_{j}\le A_{2}(\varepsilon)
    \end{displaymath}
  \end{thm}
  \begin{proof}
Indeed, for an ample Cartier divisor $H$
    \begin{displaymath}
\sum b_{j} \le \sum b_{j}B_{j}H \le -KH
    \end{displaymath}
and on each surface $X$ as above we can find $H$ with a bounded $-KH$.
  \end{proof}

  \begin{say}
The following theorem, which applies only to the case when $B$ is empty
and $-K$ is ample, but which in this case is stronger than ours,
belongs to V.V.Nikulin.
  \end{say}

  \begin{thm}[Nikulin \cite{Nikulin89a} 4.7.2]
Fix $N>0$. Then the class of surfaces with ample $-K$,
log terminal $K$ and singularities of multiplicities $\le N$ is
bounded.
  \end{thm}

\section[Boundedness for surfaces with big and nef ${K+B}$
]{Boundedness for surfaces with big and nef $K+B$}
\label{Boundedness for surfaces with big and nef K+B}

  \begin{lem}
\label{decrease one coefficient}
Assume that on a surface $X$, $K+B=K+b_{0}B_{0}+\sum_{j>0} b_{j}B_{j}$
is big and MR log canonical.
 Then one of the following is true:
    \begin{enumerate}
      \item $K+\sum_{j>0}b_{j}B_{j}$ is big;
      \item there exists $0\le b'_{0}<b_{0}$ such that
            $K+xB_{0}+\sum_{j>0}b_{j}B_{j}$ is big iff $x>b'_{0}$,
            and there exists a morphism
            $f:X\to X'$ such that
            the $\mathbb R$-divisor
            $D=f_{*}(K+b'_{0}B_{0}+\sum_{j>0}b_{j}B_{j})$
            on $X'$ is nef but not numerically trivial, and $D^{2}=0$.
            Moreover, if  $b_{0}'>0$ then the
            linear system $|ND|$, $N\gg0$ and divisible,  defines a
            $\mathbb P^{1}$-fibration $\pi:X'\to C$ to a nonsingular
            curve $C$;
     \item there exists $0\le b'_{0}<b_{0}$ such that
            $K+xB_{0}+\sum_{j>0}b_{j}B_{j}$ is big iff $x>b'_{0}$,
            and there exists a morphism
            $f:X\to X'$ such that the $\mathbb R$-divisor
            $D=f_{*}(K+b'_{0}B_{0}+\sum_{j>0}b_{j}B_{j})$
            is numerically trivial.
     \end{enumerate}
  \end{lem}
  \begin{proof}
It is well known that a divisor on a projective variety is
big if and only if it is a
sum of an ample and an effective
divisors. Hence, the property of being big is an
open property. Similarly,
the property of $(X,K+B)$ not to have a minimal model with the image
of $K+B$ nef is also an open property,
modulo the Cone and Contraction theorems of
Minimal Model Program. Indeed, it is equivalent to existence of a
covering family $\{C_{t}\}$ with $(K+B)C_{t}<0$.
Therefore, in the cases (2) and (3)
$K+b'_{0}B_{0}+\sum_{j>0}b_{j}B_{j}$ has a minimal model $X'$
with nef, but not big $D$.
If $f_{*}(B_{0})D<0$ then the Minimal Model Program applied to
$K+D-\delta B_{0}$, $0<\delta\ll1$, gives a $\mathbb P^{1}$-fibration.
In the case (2) with $b_{0}'=0$,
if the characteristic of the ground field is 0,
then some multiple of $D$ is base point free and defines an elliptic
fibration. We won't need this fact, however. In positive characteristic the
same is true but the proof requires using the classification theory,
which we would like to avoid.
  \end{proof}

  \begin{say}
The constants $A_{1}(\varepsilon)$, $A_{2}(\varepsilon)$ in the
following theorem were defined in
{}~\ref{thm A1},~\ref{thm A2 for general}.
  \end{say}

  \begin{thm}
\label{throw out almost all coefficients}
Assume that on a surface $X$, $K+B=K+\sum_{j\in J} b_{j}B_{j}$
is big and
MR $\varepsilon$-log terminal.
Then there exists a subset
of indices $J'=J'_{1}\cup J'_{2} \subset J$
such that $K+B=K+\sum_{j\in J'} b_{j}B_{j}$
is big and
    \begin{displaymath}
|J'_{1}| \le A_{1}(\varepsilon) +1 \text{, }
\sum_{j\in J'_{2}}b_{j} \le  A_{2}(\varepsilon)
    \end{displaymath}
  \end{thm}
  \begin{proof}
Decrease the coefficient $b_{0}$. We have 3
cases as in~\ref{decrease one coefficient}.

\setcounter{case}{0}
    \begin{case}
Pick another coefficient and continue.
    \end{case}

\setcounter{case}{2}
    \begin{case}
Find the maximal resolution $X''$ of $X'$ dominated by $X$:
      \begin{displaymath}
X \overset g\to X'' \overset h \to X'
      \end{displaymath}
such that
      \begin{displaymath}
K_{X''}+B''= h^{*}(f_{*}(K+b'_{0}B_{0}+\sum_{j>0}b_{j}B_{j}))
      \end{displaymath}
has all
nonnegative coefficients. Let $J'\subset J$ be the set of indices for
the nonnegative $b_{j}''$ in the latter divisor, including zeros. Then
$K_{X}+\sum_{j\in J'}b_{j}B_{j}$ is big since it contains
      \begin{displaymath}
c(K+B)+(1-c)f^{*}(f_{*}(K+b'_{0}B_{0}+\sum_{j>0}b_{j}B_{j}))\equiv
 \text{ big + nef}
      \end{displaymath}
for $0<c\ll1$. Part of the divisors $B_{j}$, $j\in J'$, lie in the
exceptional set of $X'' \to X'$ so their number is bounded by
$\rk\Pic(X'')-1$,
which is less than or equal to $A_{1}(\varepsilon)-1$ by~\ref{thm A1}.
 The sum of the coefficients for the others is bounded by
 $A_{2}(\varepsilon)$ by~\ref{thm A2 for general}. Also, we should not
 forget the divisor
 $B_{0}$ itself.
    \end{case}

\setcounter{case}{1}
    \begin{case}
Let $b'_{0}$ be as in
{}~\ref{decrease one coefficient}, (2).
To begin the proof, let us first assume that
the linear system $|ND|$, $N\gg0$ is base point free and defines a morphism
with connected fibers
$\pi\circ f:X\to C$ to a nonsingular curve $C$.
For every fiber
$F_{k}=\sum F_{kl}$ of $\pi\circ f$ consider the maximal rational number
$f_{k}$ such that
      \begin{displaymath}
K+b'_{0}B_{0}- \sum f_{k}F_{k}
      \end{displaymath}
has at least one nonnegative coefficient in $F_{kl}$, say, $F_{k0}$.
To present the    idea of the proof more clearly, assume than we can contract
for every $k$ all components of $F_{k}$ other than $F_{k0}$
to obtain a morphism
      \begin{displaymath}
f''':X\to X'''
      \end{displaymath}
Note that the coefficients of components of
      \begin{displaymath}
K_{X'''}+B'''= f'''_{*}(K+b'_{0}B_{0} -\sum f_{k}F_{k})
      \end{displaymath}
are all nonnegative.
      \begin{subcase}$K_{X'''}+B'''$ is nef but not numerically trivial.

Then $K_{X}+B'''+ (b_{0}-b'_{0})B_{0}$ is big,
so we can use $J'$ corresponding to $B'''$
plus $B_{0}$ itself, i.e.\  to all $B_{j}$ not
in fibers of $\pi\circ f$. Evidently,
        \begin{displaymath}
\sum_{j\in J'-0}b_{j} \le 2 < A_{2}(\varepsilon)
        \end{displaymath}
      \end{subcase}

      \begin{subcase} The opposite to the subcase 1.

At the same time $K_{X'''}+B'''+\sum b_{k0}B_{k0}$ is nef but not big.
Hence, there
exist several $F_{k0}$, say, $F_{00},F_{10}\ldots{}F_{s0}$ such that
        \begin{displaymath}
K_{X'''}+B'''+b'_{00}F_{00}+ \sum_{i=1}^{s} b_{i0}B_{i0}
        \end{displaymath}
is numerically trivial.
Find the maximal partial resolution $X''$ of $X'''$ dominated by $X$:
      \begin{displaymath}
X \overset g\to X'' \overset h \to X'''
      \end{displaymath}
such that
      \begin{displaymath}
K_{X''}+B''= h^{*}(K_{X'''}+B''' +b'_{00}F_{00}+\sum_{i=1}^{s}b_{i0}B_{i0}
      \end{displaymath}
has all nonnegative coefficients.
Then choose $J'$ as in the case 3. Again, we can divide $J'$ into
two parts, $J'_{1}$ and $J'_{2}$ and count indices as above. We
shouldn't forget one
more divisor, $F_{00}$, to arrive at the final estimate.
      \end{subcase}

In this proof we assumed the existence of a fibration $\pi\circ f:X\to
C$. If there is a component $B_{k}$ of $b_{0}'B_{0}+\sum_{j>0}b_{j}B_{j}$
that intersects $D$ positively (i.e.\ ``horizontal'')
then applying the Cone and the
Contraction theorems to $K+b_{0}'B_{0}+\sum_{j>0}b_{j}B_{j}-\epsilon
B_{k}$ for $0<\epsilon\ll1$ several times we obtain such a $(\mathbb
P^{1})$-fibration.
The case when no such component $B_{k}$ exists corresponds to an elliptic
fibration. We could still prove that a multiple of $D$ is
base point free, i.e.\  the abundance theorem,
(cf.~\cite{TsunodaMiyanishi83}), but this would involve some
classification theory. Instead, we can use the same argument as above
with not actual fibers but fibers in the numerical sense. Consider
a connected component of $\supp(b_{0}'B_{0}+\sum_{j>0}b_{j}B_{j})$.
Then by the Hodge Index theorem the corresponding graph is elliptic or
parabolic.
A parabolic graph corresponds to a fiber, an elliptic graph --
to a part of a fiber, the
case when several curves of the fiber have coefficients 0.
Then argue as above.

As for the existence of a contraction of components of $F_{k}$ other
than $F_{k0}$,
note that all we need in this
proof is the contraction $g:X\to X''$ and it can be
constructed without obtaining $X'''$  first. Now, $g:X\to X''$ exists
because the corresponding
graph is log
terminal, cf.~ section~\ref{MMP in dimension 2}.

    \end{case}
\renewcommand{\qed}{}
  \end{proof}

  \begin{cor}
\label{cor to throw out almost all coefficients}
    \begin{displaymath}
|J'| \le A_{1}(\varepsilon)+A_{2}(\varepsilon)/\min(b_{j}) +1
    \end{displaymath}
  \end{cor}

  \begin{thm}
\label{can decrease all by delta}
Assume that on a surface $X$, $K+B=K+\sum b_{j}B_{j}$
is big and MR log canonical.
Further assume that $b_{j}$ belong to a D.C.C. set $\mathcal C$. Then there
exists a constant
$\delta(\mathcal C)>0$ depending only on the set $\mathcal C$ such that
$K+\sum (b_{j}-\delta)B_{j}$ is big.
  \end{thm}
  \begin{proof}
Assume the opposite. Then, similarly to the proof of
{}~\ref{decrease one coefficient}, there exists an
infinite sequence $\delta_{n}\to0$ and a sequence of surfaces $X'_{n}$
such that the direct image of
$K+\sum (b_{j}-\delta_{n})B_{j}$,
which we will denote by $K+B'_{n}$,
 is numerically trivial on a
general fiber of a $\mathbb P^{1}$-fibration or simply is numerically
trivial. In the former case for some nonnegative integers $k_{j}$ we get
    \begin{displaymath}
\sum k_{j}(b_{j}-\delta_{n})=2
    \end{displaymath}
which is trivially impossible. In the latter case we get a
contradiction  as follows. If $B'_{n}$ is
not empty then $K+B'_{n}-\alpha B_{1}$ is not nef and the corresponding
extremal ray does not contain $B_{1}$. Therefore, after several
contractions we can assume that $\rk\Pic X'_{n}=1$. Then the number
of irreducible components in $B'_{n}$ is less than $3/\min(\mathcal C)+1$
and the coefficients of $B'_{n}$ give an infinite
strictly increasing sequence of chains $\{b_{n,k}\}$. This is impossible by
\cite{Alexeev93a}, 5.3. Analyzing the proof of 5.3 shows that the
condition for $K+B'_{n}$ to be log canonical can be weakened to MR log
canonical.
  \end{proof}

  \begin{thm}
\label{main diagram}
Fix $C>0$ and a D.C.C. set $\mathcal C$.
Then there exist a bounded class of
surfaces with divisors $(Z,D)$ such that for every surface $X$ with
$K+B$ nef, big, MR log canonical, with $b_{j}\in \mathcal C$ and
$(K+B)^{2}\le C$
there exists a diagram
    \begin{displaymath}
    \begin{CD}
      Y @>g>> Z\\
      @VfVV\\
      X
    \end{CD}
    \end{displaymath}
in which
    \begin{enumerate}
      \item $Y$ is the minimal desingularization of $X$,
      \item defining $K_{Y}+B^{Y}=f^{*}(K+B)$,
            $D=g(\supp B^{Y}\cup \Exc(f))$ where $\Exc(f)$ is the union of
            exceptional divisors of $f$.
    \end{enumerate}
  \end{thm}
  \begin{proof}

Start with a divisor $K_{Y}+f^{-1}B+\sum F_{i}$ on $Y$.
By~\ref{can decrease all by delta} we can decrease the coefficients of
$B^{Y}$ to obtain a divisor $B'$ with the following properties:
    \begin{enumerate}
      \item $K_{Y}+B'$ is big,
      \item coefficients of $B'$ belong to a {\em finite} set
            $\mathcal C'$ of {\em rational} numbers that depends only on
            the original set $\mathcal C$,
      \item all the coefficients of $B'$ are less than
      $1-\varepsilon<1$, i.e.\  $K_{Y}+B'$ is MR $\varepsilon$-log
      terminal. Here $\varepsilon$ depends again only on $\mathcal C$.
    \end{enumerate}

Next use~\ref{throw out almost all coefficients} to obtain a big
divisor $K_{Y}+B''$. The coefficients of $B''$ belong to the same set
$\mathcal C'$ and, moreover, the number of components of $B''$ is bounded
in terms of $\mathcal C$.

Apply the log Minimal Model Program in a slightly more general than
usual form
described in section~\ref{MMP in dimension 2} (because $K_{Y}+B''$ is
not necessarily log canonical)
to obtain a log canonical model $g:Y\to Z$
of $K_{Y}+B''$. Then $H=g(K_{Y}+B'')$ is ample and gives a $\mathbb
Q$-polarization on $Z$.
{}From the classification of surface log terminal singularities
(cf. \cite{Alexeev93a} 3.3)
it follows
that $Z$ can have only finitely many types of singularities.
Indeed, for each of these singularities number of exceptional curves
on the minimal
resolution is bounded, and all log discrepancies are
$\ge\varepsilon>0$,
so the corresponding weights are less than $2/\varepsilon$.
Therefore
there exists an integer $N(\mathcal C)$ such that $NH$ is a polarization
of $Z$. Note that by construction $H\underset{c}\le K+B$,
hence
    \begin{displaymath}
H^{2}\le (K+B)^{2} \le C
    \end{displaymath}
by~\ref{easy inequalities}.
Also, $HK\le H^{2}$ and
    \begin{displaymath}
p_{a}(NH)=1/2(N^{2}H^{2}+NHK)\ge0.
    \end{displaymath}
Hence, there are only finitely many possible values for $H^{2}$ and $HK$.
Now we use~\ref{main method for proving boundedness} to
conclude that the class of surfaces $\{Z\}$ is bounded. To prove that
the class $\{(Z,D)\}$ is bounded we also
have to show that  $DH=g(\supp B\cup \Exc(f))H$ is bounded.
This is clear for the components of $g(\supp B'')$. So let us
consider the rest of $g(\supp B\cup \Exc(f))$.
First, consider  divisors $D_{k}$ on $Y$
with $g(D_{k})\ne\text{point}$ which are not
$f$-exceptional $(-2)$-curves.
Their coefficients $d_{k}$ in $K_{Y}+B^{Y}$ are bounded from below by
$\min(1/3,\mathcal C)$. We have  $g^{*}(H)\le K_{Y}+{B^{Y}}$ by
{}~\ref{inequalities}, so we
get the desired bound from
    \begin{eqnarray*}
(K_{Y}+B^{Y}-g^{*}(H))g^{*}(H) \le
(K_{Y}+B^{Y})g^{*}(H) \le
(K_{Y}+B^{Y})^{2} \le C
    \end{eqnarray*}

 For $f$-exceptional (-2)-curves one has
    \begin{displaymath}
E=g^{*}(H)+1/2\sum(Hg(D_{k}))D_{k}\underset{c}{\le}K+B^{Y}
    \end{displaymath}
by~\ref{inequalities} because we add to $g^{*}(H)$
only $f$-exceptional divisors and
    \begin{displaymath}
E^{2}\ge H^{2}+1/2\sum(Hg(D_{k}))^{2}.
    \end{displaymath}
On the other hand,
by lemma~\ref{only surfaces}
    \begin{displaymath}
E^{2}\le (K_{Y}+B^{Y})^{2}\le C.
    \end{displaymath}
This concludes the proof.
  \end{proof}

  \begin{thm}
\label{bound for varepsilon-terminal, K big, nef}
Fix $\varepsilon>0$, a constant $C$ and a D.C.C. set $\mathcal C$.
          Consider all surfaces $X$
          with an $\mathbb R$-divisor $B=\sum b_{j}B_{j}$ such that
          $K_{X}+B$ is MR $\varepsilon$-log terminal, $K_{X}+B$ is
          big and nef, $b_{j}\in \mathcal C$ and $(K+B)^{2}\le C$.
          Then the class $\{(X,\supp B)\}$ is bounded
  \end{thm}
  \begin{proof}
Consider a diagram as in~\ref{main diagram} above.
Because the family $(Z,D)$ is bounded,
changing $Z$ and $D$ we can assume
that $Z$ is nonsingular and that the only singularities of $D$
in the exceptional set for $Y\to Z$ are nodes.
But then $Y$ is obtained from $Z$ by several
blowups at these nodes of $D$
and their number depends only on $\varepsilon$. Indeed,
since $K_{Y}+B^{Y}$ is nef, when
blowing up the point of intersection of two curves with coefficients
$b_{1}$ and $b_{2}$, the new coefficient should be $\le
b_{1}+b_{2}-1$. On the other hand, all coefficients in  $K_{Y}+B^{Y}$
are nonnegative. Since every $b_{k}\le1-\varepsilon<1$, only finitely
many such blowups can be done.

By~\ref{boundedness by blowingup}
 we conclude that the family $(Y,\supp B\cup \Exc(f))$ is bounded.
The surfaces $Y$ come with a polarization, call it $H_{Y}$ with
bounded $H_{Y}^{2}$, $H_{Y}K_{Y}$, $H_{Y}(\supp B\cup \Exc(f))$.
Since there are only finitely many configurations
for the exceptional curves $\Exc(f)$ and they are all log terminal,
hence rational, we can construct a polarization $H_{X}$ on $X$ with
bounded  $H_{Z}^{2}$, $H_{Z}K_{Z}$, $H_{Z}(\supp B)$.
At this point, we can use~\ref{main method for proving boundedness}
one more time.
  \end{proof}

\section[Descending Chain Condition]{Descending Chain Condition}

  \begin{say}
The aim of this section is to prove the following
  \end{say}

  \begin{thm}
\label{DCC for semi-log canonical}
Fix a D.C.C. set $\mathcal C$.
Consider all surfaces $X$
with an $\mathbb R$-divisor $B=\sum b_{j}B_{j}$ such that
$K_{X}+B$ is MR semi-log canonical, $K_{X}+B$ is
ample and $b_{j}\in \mathcal C$.
Then the set $\{ (K_{X}+B)^{2} \}$ is a D.C.C. set.
  \end{thm}

  \begin{say}
The theorem will be proved in several steps, under weaker and weaker
conditions.
  \end{say}

  \begin{thm}
\label{DCC for Kawamata log terminal}
 Theorem~\ref{DCC for semi-log canonical} holds for $K+B$
MR Kawamata log terminal.
  \end{thm}
  \begin{proof}
We can certainly assume that there exists a
constant $C$ such that $(K+B)^{2}\le C$. By
{}~\ref{main diagram}, we know
that there exists a diagram
    \begin{displaymath}
    \begin{CD}
      Y @>g>> Z\\
      @VfVV\\
      X
    \end{CD}
    \end{displaymath}
with $Z$ and $D=\supp B\cup \Exc(f)$ bounded.

Assume first that $Z$ and $D$ are actually fixed. Then the proof immediately
follows from~\ref{kinky method for proving DCC} and the following
theorem, since for surfaces $H^{2} \le (g^{*}g_{*}H)^{2}$.

  \begin{thm}
\label{thm main trick with sequences}
    Let $\{Y^{(s)} \,|\, s\in \mathbb N \}$ be a sequence of surfaces as above.
    Then changing $(Z,D)$ and picking a subsequence one has for every $s<t$
    \begin{equation}
\label{votono}
      g^{(t)*}(g^{(s)}_{*}(K+B^{Y^{(s)}})){\le}
      K+B^{Y^{(t)}}
    \end{equation}
  \end{thm}
  \begin{proof}
\setcounter{case}{0}
Let me outline the strategy first.
After
canceling $K$ on both sides of the inequality~\ref{votono}, we have to
compare the coefficients for all prime divisors.
We are going to look at
coefficients in $B^{Y^{(t)}}$.
Some of them come from $B_{X^{(t)}}$ and hence
belong to the set
$\mathcal C$. Since this set satisfies a descending chain condition,
we can hope that after passing to a subsequence all the
required inequalities will be satisfied.
Then, the other part of coefficients in $B^{Y^{(t)}}$ comes from the
exceptional divisors of $f^{(t)}$. {\em A priori}, there is no
information about these coefficients except that they are nonnegative
and are less than 1. But, by
lemma~\ref{inequalities} as soon as we prove the inequalities for the
first part of $B^{Y^{(t)}}$, the inequalities for the exceptional
divisors will follow automatically. Finally, the fact that
$K+B^{Y^{(s)}}$ is MR Kawamata
log terminal will mean that for every fixed $s$, after modifying  $Z$
and $D$, there are only {\em finitely many} inequalities we have to
check.

Consider $g^{(s)}:Y^{(s)}\to Z$ and the divisors
$g_{*}(K_{Y^{(s)}}+B^{Y^{(s)}})$ on $Z$. We shall denote
various prime divisors by $D_{\ldots{}}$ and the corresponding
sequences of coefficients in $K+B^{Y^{(s)}}$ by $x^{(s)}_{\ldots{}}$,
where $\ldots{}$ stands for a certain index.
For any fixed divisor $D_{\ldots{}}$, after
picking a subsequence, we can
assume that the sequence $x^{(s)}_{\ldots{}}$ is either increasing or
decreasing (not necessarily strictly) and has the limit
$c_{\ldots{}}$. Below we shall talk, for example, about a divisor
$D_{\ldots{}}$ that we get by blowing up a particular point $P$ on $Z$.
Picking a subsequence, we can assume that this point is blown up for
all the surfaces $Y^{(s)}$ or is not blown up at all
and then the analysis of the coefficients $x^{(s)}_{\ldots{}}$
is superfluous. The important fact that we shall use is that all the
coefficients $x_{\ldots{}}$ are strictly less than 1.

    \begin{step}
Changing $(Z,D)$ we can assume that the points being blown up are
nonsingular on $Z$ and are at
worst the nodes of $D$. Moreover, if we blow up a point $P$ which is better
than a node then for any divisor $D_{k}$ over $P$ in $B^{Y^{(t)}}$ the
corresponding coefficient in
$g^{(t){*}}(g_{*}^{(s)}(K+B^{Y^{(s)}}))$
is negative, so we have the inequality~\ref{votono}.

Therefore, below we consider only the situation of two normally
crossing divisors $D_{1}$ and $D_{2}$ and the corresponding numbers
$c_{1},c_{2}\le1.$
    \end{step}

    \begin{step}
Points $P$ with $c_{1}=c_{2}=1$.

Blow up $P$ to get $D_{3}$. If $c_{3}<1$, change $Z$ to this new blown
up surface and go on to the next step. If $c_{3}=1$ then for any divisor
$D_{k}$ over $P$ we have $c_{k}=1$ by a very simple computation taking
into account that all $K+B^{Y^{(s)}}$ are nef.
Since the coefficients $\xso$ and
$\xst$ are both less than 1, for the inequality~\ref{votono} to be
true for $s=1$
over $P$ we need to check only finitely many divisors and the
inequalities are satisfied after picking a subsequence
because $\underset{s\to\infty}{\lim} \xsk=c_{k}=1$. Then do the same for
$s=2$ and so on. For the new sequence the inequalities~\ref{votono} are
satisfied for all divisors over $P$.
    \end{step}

    \begin{step}
Points $P$ with $c_{1}=1$, $c_{2}<1.$

Let us change notation slightly in this step. Denote $D_{2;1}=D_{2}$,
and then
$D_{2;k+1}$ will be the exceptional divisor of a blowup at $D_{1}\cap
D_{2;k}$, $k=1,2,3\ldots{}$ Then denote by $D_{3;k}$ the exceptional
divisor of a blowup at the point $D_{2;k}\cap D_{2;k+1}$, by $D_{4;k}$
the exceptional divisor of a blowup at $D_{2;k}\cap D_{3;k}$ (or at
$D_{3;k}\cap D_{2;k+1}$) and so on.

The claims below are elementary and follow immediately from the fact
that $K+B^{Y^{(s)}}$ are nef.

      \begin{claim}
        In this way we obtain {\em finitely\/}
        many sequences $\{ \xspk \,|\,
        k\in \mathbb N     \}$ that are the only possible nonnegative
        coefficients in $B^{Y^{(s)}}$.
      \end{claim}

Note that $c_{2;1}\ge c_{2;2}\ge \ldots{}$ and denote by
$c=\underset{k\to\infty}{\lim}c_{2;k}$.

      \begin{claim}
        There exist naturally defined positive integers $n(p)$ such that
        \begin{displaymath}
        \underset{k\to\infty}{\lim}c_{p;k}= 1-n(p)(1-c)
        \end{displaymath}
      \end{claim}

      \begin{claim}
\label{convex}
If for some $k,s$  $x^{(s)}_{2;k+2}< x^{(s)}_{2;k+1}>d$ then for the same $k,s$
one also has $\xspk>1-n(p)(1-d)$
       \end{claim}

\bigskip
Now consider two cases.

       \begin{case}
There exists $k_{0}$ such that for infinitely many $s$
         \begin{displaymath}
           \xstkz\le c
         \end{displaymath}

Pick a subsequence of $s$ with this property. Change $Z$ by its
$(k_{0}-1)$-blowup. Then again for the inequality~\ref{votono} to be
true for $s=1$ we have to check only finitely many divisors and we get the
corresponding inequalities because $c_{p;k}\ge 1-n(p)(1-c)$
by~\ref{convex} with $d=c$. Then proceed
with $s=2$ and so on.
       \end{case}

       \begin{case}
For any $k$ there exist only finitely many $s$ such that
         \begin{displaymath}
           \xstk\le c
         \end{displaymath}
Then~\ref{convex} with $d=c$ implies that
for any fixed $p,k$ there exist only finitely many $s$ such that
         \begin{displaymath}
           \xspk\le 1-n(p)(1-c)
         \end{displaymath}

Define a set
         \begin{displaymath}
\mathcal A_{2}= \{ \xstk \,|\, \xstk>c; \, k,s\in \mathbb N \}
         \end{displaymath}
Fix $k_{0}$ such that $c_{2;k_{0}}< \min\mathcal A_{2}\cap \mathcal C$, which
is possible because $\mathcal C$ is a D.C.C. set. Then pick a subsequence
of $s$ with the following property:

         \begin{equation}
x\in \mathcal A_{2}, \, x< \xstkz \text{ implies } x\notin \mathcal C
         \end{equation}

Introducing
         \begin{displaymath}
\mathcal A_{p}= \{ \xspk \,|\, \xspk>1-n(p)(1-c); \, k,s\in \mathbb N \}
         \end{displaymath}
we can also assume that
         \begin{equation}
x\in \mathcal A_{p}, \, x< 1-n(p)(1-\xstkz) \text{ implies } x\notin \mathcal C
         \end{equation}

Now, again, as above, we change $Z$ by its $(k_{0}-1)$-blowup and arrange
the inequality~\ref{votono} to be true for $s=1$ picking a
subsequence. The important thing here is that the relevant inequalities
may fail only for $x\notin \mathcal C$ which is OK
by lemma~\ref{inequalities}.

Then proceed with the same procedure for $s=2$ and so on.
       \end{case}
Checking for a vicious circle, I emphasize that in this step we change
$Z$ for each point $P$ only once.
    \end{step}

    \begin{step}
Points with $c_{1}<1$, $c_{2}<1$.

Change $Z$ by blowing it up several times so that for any future
blowup all divisors should have negative coefficients in $B^{Y^{(s)}}$.
    \end{step}

At this point we achieved that the inequality~\ref{votono} is satisfied
for any divisor which is exceptional for $g^{(s)}$. Finally, let us
settle it for the divisors $D_{k}$ on $Z$ itself.

    \begin{step}
$D_{1}$ is a divisor on $Z$.

If $\{ \xso \}$ is an increasing sequence, we are done. Otherwise,
omitting finitely many $s$, we have $\xso< \min \{\xso\} \cap \mathcal C$, so
the coefficients $\xso$ correspond to exceptional divisors and
we are OK by lemma~\ref{inequalities}.
    \end{step}

So far, we assumed that $(Z,D)$ is fixed.
To prove the general case we use~\ref{kinky method for proving DCC}. Pick
an arbitrary sequence of surfaces satisfying the initial conditions.
Since  by~\ref{main diagram}
the class $(Z,D)$ is bounded, there are only
finitely many numerical possibilities for $Z$, $D$ and the subsequent
blowups so, after picking a subsequence,
 we can argue in the same way as we did for the fixed $(Z,D)$.
\renewcommand{\qed}{}
\qed
  \end{proof}

  \begin{rem}
Several important ideas of
the proof given here are taken directly
 from the Xiao's proof of
Proposition~5, \cite{Xiao91}.
Unfortunately, I do not understand the arguments of
\cite{Xiao91} completely and therefore cannot say precisely
how close I follow them. However, it is
clear that
the Xiao's proof may apply only to a set $\mathcal C$
with a single limit point at 1.
  \end{rem}

  \begin{thm}
\label{DCC for log canonical}
Theorem~\ref{DCC for semi-log canonical} holds for $K+B$ MR log
canonical.
  \end{thm}
  \begin{proof}
Let $\{X^{(s)},K+B^{(s)}\}$ be a sequence with strictly decreasing
$(K+B)^{2}$. For each $X^{(s)}$ consider the maximal log crepant
(i.e.\ all new log codiscrepancies equal~1)
partial resolution $f^{(s)}$ which is dominated by the minimal
desingularization. Let $E_{i}$ be the divisors with coefficients~1 in
$f^{*}(K+B)$. Then for the appropriate choice of
$\epsilon_{i}\to0$ the
divisors $f^{*}(K+B)-\sum \epsilon_{i}E_{i}$ will be ample and
MR Kawamata log
terminal, the coefficients will belong to a {\em new}, but again a
 D.C.C., set,
and the sequence of the squares will still be strictly decreasing.
So we get a contradiction by~\ref{kinky method for proving DCC}
and \ref{DCC for Kawamata log terminal}.
  \end{proof}

{\em End of the proof of theorem~\ref{DCC for semi-log canonical}.}

By~\ref{reduce semi-log canonical to log canonical},
  \begin{displaymath}
(K+B)^{2}= \sum (K_{X_{m}}+B_{m})^{2}
  \end{displaymath}
and all $K_{X_{m}}+B_{m}$ are log canonical. By
{}~\ref{DCC for log canonical} every summand in this formula belongs
to a D.C.C. set, and so does the sum.
  \end{proof}

\section[Boundedness for the constant $(K+B)^{2}$
]{Boundedness for the constant $(K+B)^{2}$}

  \begin{say}
In this section we shall prove the following theorem. Again, it will
be done first for $K+B$ MR Kawamata log terminal and MR log canonical.
  \end{say}

  \begin{thm}
\label{K2=const for semi-log canonical}
Fix a constant $C$ and a D.C.C. set $\mathcal C$.
Consider all surfaces $X$
with an $\mathbb R$-divisor $B=\sum b_{j}B_{j}$ such that
$K_{X}+B$ is MR semi-log canonical, $K_{X}+B$ is
ample, $b_{j}\in \mathcal C$ and $(K+B)^{2}= C$.
Then the class $\{(X,\sum b_{j}B_{j})\}$ is bounded.
  \end{thm}

  \begin{thm}
Theorem~\ref{K2=const for semi-log canonical} holds for $K+B$ MR
Kawamata log terminal. Moreover, in this case $K+B$ may be taken to
be big and nef instead of ample.
  \end{thm}
  \begin{proof}
Again, first assume that in the diagram~\ref{main diagram} $(Z,D)$ is
fixed. Then, for the general case we use~\ref{kinky method for proving
boundedness} since $(Z,D)$ moves in a bounded family.

By lemma~\ref{inequalities} the inequalities~\ref{votono} become
equalities only when
    \begin{displaymath}
      g^{(s)*}(K+\overline B)=\kbys \, \text{ for every }\, s \,
      \text{ with a fixed }\, \overline B
    \end{displaymath}
Since the log canonical model of $\kbys$
 is the image of the linear system $|N(\kbys)|$, $N\gg0$,
all surfaces in the sequence have the same log canonical model.
This means that the class $\{(Z,\overline{ B})\}$ is bounded and since
all coefficients in $B^{Y^{(s)}}$ are less than~1, that $\{(Y,B)\}$ is
bounded.
Since
    \begin{displaymath}
K_{Y^{(s)}}+B^{Y^{(s)}}= f^{*}(K_{X^{(s)}}+B_{X^{(s)}}),
    \end{displaymath}
we get a diagram $Y^{(s)}\to X^{(s)} \to Z$
and the class $\{(X,B)\}$ is bounded by the ``sandwich'' principle.
  \end{proof}

  \begin{thm}
\label{K2=const for log canonical}
Theorem~\ref{K2=const for semi-log canonical} holds for $K+B$ MR log
canonical.
  \end{thm}
  \begin{proof}
As in the previous proof,
by~\ref{main diagram} we first assume existence of the diagram
    \begin{displaymath}
    \begin{CD}
      Y @>g>> Z\\
      @VfVV\\
      X
    \end{CD}
    \end{displaymath}
with $Z$ and $D=\supp B\cup \Exc(f)$ fixed.
Pick an infinite sequence of surfaces $\{ X^{(s)} \}$ and for
$K_{Y^{(s)}}+B^{Y^{(s)}}=f^{(s)*}(K_{X^{(s)}}+B^{(s)})$
write
    \begin{displaymath}
(K_{Y^{(s)}}+B^{Y^{(s)}})^{2}= Base^{(s)} + Defect^{(s)}
    \end{displaymath}
where
  \begin{eqnarray*}
  Base^{(s)}= (g^{(s)}_{*}(K_{Y^{(s)}}+B^{Y^{(s)}}))^{2}
\\Defect^{(s)}= (K_{Y^{(s)}}+B^{Y^{(s)}})^{2}-
(g^{(s)}_{*}(K_{Y^{(s)}}+B^{Y^{(s)}}))^{2}
  \end{eqnarray*}
Note that $Defect^{(s)}\le0$ for all $s$.

Now if $K+B$ were MR Kawamata log terminal then from
{}~\ref{thm main trick with sequences} we could
conclude that, after picking a subsequence, $\{ Base^{(s)} \}$ is a
nondecreasing sequence and
that $Defect^{(s)}\to0$. But because every MR log
canonical $K+B$ can be approximated by MR Kawamata log terminal
divisors as in the proof of
{}~\ref{DCC for log canonical}, we can conclude just the same in the present
situation. But then $Base^{(s)} + Defect^{(s)}=C$ implies that
$Defect^{(s)}=0$ and $Base^{(s)}=C$,
i.e.\
    \begin{displaymath}
      g^{(s)*}(K+\overline B)=\kbys \, \text{ for every }\, s \,
      \text{ with a fixed }\, \overline B
    \end{displaymath}
Since $\kbxs$
 is the image of the linear system $|N(\kbys)|$, $N\gg0$,
we get $(X^{(s)},B)=(Z,\overline{ B})$, so
all surfaces in the sequence are actually isomorphic to each other.
In the general situation $(Z,D)$ is not fixed but the class $(Z,D)$ is
bounded, so we get that the class $(X,B)$ is bounded.
    \end{proof}

{\em End of the proof of~\ref{K2=const for semi-log canonical}. }
    \begin{proof}
Let $X^{\nu}=\cup X_{m}$ be the normalization of $X$. Then as in
{}~\ref{reduce semi-log canonical to log canonical}
      \begin{displaymath}
(K+B)^{2}= \sum_{m} (K_{X_{m}}+B_{m})^{2},
      \end{displaymath}
the coefficients of $B_{m}$ belong to $\mathcal C\cup \{1\}$ and
$K_{X_{m}}+B_{m}$ are log canonical.
Applying~\ref{DCC for log canonical} and~\ref{K2=const for log
canonical} we see that there are only finitely many
possibilities for $(K_{B_{m}})^{2}$ and that all
$(X,\supp B_{k})$ belong to a bounded class. Now the class $\{X\}$ is
bounded by the Conductor Principle \cite{Kollar85} 2.3.5. It
immediately follows that in fact  $\{(X,\supp B)\}$ is bounded.
    \end{proof}

\section[On the log Minimal Model Program for surfaces,
in arbitrary characteristic
]{On log MMP for surfaces}
\label{MMP in dimension 2}

    \begin{say}
The Log Minimal Model Program in dimension two is undoubtedly much
easier than in higher dimensions. There are two main circumstances that
distinguish the surface case.

The first difference is that
the Log Minimal
Model Program in dimension two is characteristic free, see
\cite{Mori82} for nonsingular surfaces and \cite{TsunodaMiyanishi83}
for open surfaces. The Cone Theorem follows from the nonsingular case
by a simple argument as in \cite{TsunodaMiyanishi83}. The proofs of
the Contraction Theorem and of the Log Abundance Theorem in higher
dimension use the Kodaira vanishing theorem which fails in positive
characteristic. However, in dimension two they can be proved
using the contractibility conditions of Artin \cite{Artin62} because
the corresponding configurations of curves are rational. In fact, in
characteristic $p>0$ contractibility conditions are even weaker than
in characteristic 0.

The second difference is the existence of a unique minimal resolution of
singularities.

The statements in this section and their proofs are elementary.
I provide
them here because I need them for this paper and they are slightly
more general than those appearing elsewhere.

	\end{say}

    \begin{thm}[The Cone Theorem]
Let $X$  be a normal surface and  $K+B=K_{X}+\sum b_{j}B_{j}$ be an
$\mathbb R$-Cartier divisor with $b_{j}\ge0$. Let $A$ be an ample  divisor on
$X$. Then for any $\epsilon>0$ the Mori-Kleiman cone of effective
curves $\overline{NE}(X)$ in $NS(X)\otimes\mathbb R$ can be written as
      \begin{displaymath}
\overline{NE}(X)= \overline{NE}_{K+B+\epsilon A}(X) +\sum R_{k}
      \end{displaymath}
where, as  usual, the first part consists of cycles that have positive
intersection with $K+B+\epsilon A$ and  $R_{k}$ are finitely many
extremal rays. Each of the extremal rays is generated by an effective curve.
    \end{thm}
    \begin{proof}
It suffices to prove this theorem for the minimal desingularization
$Y$ of
$X$ with $K_{Y}+B^{Y}=f^{*}(K+B)$ instead of $K+B$ since every effective
curve
on $X$ is an image of an effective curve on $Y$,
and so $\overline{NE}(X)$ is the image of $\overline{NE}(Y)$ under a
linear map $f_{*}:N_{1}(Y)\to N_{1}(X)$.
 Then for a
nonsingular surface the statement easily follows from
the usual Cone Theorem
\cite{Mori82} by the same argument as  in \cite{TsunodaMiyanishi83}, 2.5.
The reason for this  is,   certainly, that for any curve $C\ne B_{j}$
one has $C(K+B)\ge CK$.
The set of
extremal rays   for $K+B+\epsilon A$ is a subset of extremal rays for
$ K+\epsilon A$, possibly enlarged by some of the  curves $B_{j}$.
    \end{proof}

    \begin{thm}[Contraction Theorem]
Let $X$ be a projective surface with MR log canonical
$K+B=K+ \sum  b_j B_j$.  Let $R$ be an
extremal ray for $K+B$. Then there exists a nontrivial projective
morphism $\phi_{R}:X\to Z$ such that $\phi_{R*}(\mathcal O_{X})=\mathcal
O_{Z}$ and $\phi(C)=pt$ iff the class of $C$ belongs to $R$, and
$K_{Z}$ is also log canonical. If $K_{X}$ is Kawamata log terminal
then so is $K_{Z}$.
    \end{thm}
    \begin{proof}
By the previous theorem $R$ is generated by an irreducible curve $C$.
If $C^{2}>0$ then  as in  \cite{Mori82} 2.5, $f^{*}(C)$ on the minimal
resolution  (defined, according to Mumford, for every Weil divisor
$C$) is in the interior of $\overline{NE}(Y)$, so $C$ is in the
interior of $\overline{NE}(X)$. This is possible only if $\rho(X)=1$
and then
$\phi_{R}$ maps the whole $X$ to a    point.

If $C^{2}=0$ then
the graph corresponding to $\supp(f^{*}(C))$ is parabolic and it
is rational because $KC<0$.
By the Riemann--Roch formula $h^{0}(nf^{*}(C))$ grows linearly in   $n$,
so for some  $N\gg0$ the linear system $|Nf^{*}(C)|$ is base point free  and
defines a projective morphism to a curve. Similarly  to \cite{Mori82},
2.5.1, this is   possible only if   $\rho(X)=2$.

Now let us assume that $C^{2}<0$. Consider the minimal resolution
$f:Y\to X$ and
let $F_{i}$ be  the maximal number of exceptional  curves of $f$ such
that $C\cup_{i} F_{i}$ is connected. Then, since $(K+B)C<0$ and
the quadratic form $|F_{i}F_{i'}|$ is negative definite, all the
log    discrepancies of the graph $C\cup_{i} F_{i}$ are strictly
positive, so this  graph is log terminal. Looking at the list of the
minimal log
terminal graphs (for example, in \cite{Alexeev92}) one easily observes that
they
are all rational in the sense of Artin \cite{Artin62}.
So   are all the nonminimal
log terminal elliptic graphs obtained from them by simple
blowups.
Hence, by
\cite{Artin62} the configuration of  curves $C\cup_{i}F_{i}$ can be
contracted to a normal point on a projective surface $Z$. By normality
of $X$, this defines a morphism $\phi_{R}:X\to Z$ satisfying the
required properties.
    \end{proof}

    \begin{thm}[Easy Log Abundance Theorem]
Let $X$ be a projective surface with MR Kawamata log terminal
$K+B=K+\sum b_{j}B_{j}$.
Assume
that $K+B$ is   big and nef. Then for some  $N\gg0$ the linear system
$|N(K+B)|$ is base point free and defines a birational projective
morphism $f:X\to Z$ with $N(K+B)=f^{*}(H)$ for some   ample divisor
$H$ on   $Z$. $K_{Z}$ is also Kawamata log terminal.
    \end{thm}
    \begin{proof}
We have $(K+B)^{2}>0$.
If $K+B$ is not ample then by the Nakai-Moishezon criterion of ampleness
there exists a  curve $C$ such that $(K+B)C=0$.
By the Hodge Index theorem, $C^{2}<0$.
Then, as  in the
previous theorem, $C$ can be     contracted to a  normal point on a
projective surface $Z_{1}$, $K+B=g^{*}(K_{Z_{1}}+g_{*}(B))$ and
the  latter divisor is again MR Kawamata log terminal.
Now the statement follows
by induction  on rank of the Picard group of $Z_{1}$.
    \end{proof}

    \begin{rem}
It is elementary to generalize the  above three theorems to a relative
situation of an arbitrary
projective  morphism $X\to S$ to a variety $S$.
    \end{rem}

\section{Concluding remarks, generalizations, open questions}
\label{Concluding remarks, generalizations, open questions}

  \begin{say}
Following the logics of the theorems presented here one should expect
that there exist infinitely many types of surfaces with ample
canonical divisor $K$, $K^{2}\le   C$
 and {\em empty\/} $B$ if we allow singularities worse
than $\varepsilon$-log terminal,
for example arbitrary quotient singularities. It would be interesting
to see the  examples.
  \end{say}

  \begin{say}[Effectiveness]
Note that it is quite straightforward to obtain effective formulas for
$A_{1}$ in~\ref{thm A1} and $S_{1}$ in~\ref{minimal elliptic graph}.
However, methods of this paper do not provide
effective estimates for functions $A_{2}$ in~\ref{thm A2 for
general} and $c(\mathcal C)$ in~\ref{my thm on auts}.
  \end{say}

  \begin{say}
I would go as far as to conjecture that the conditions for the
boundedness formulated in~\ref{proved here}
 are the most natural ones and that direct
analogs of all four our main theorems
should be true in any dimension.
    \end{say}

    \begin{say}
In the case of $K$ negative, i.e.\ Fano varieties, there are several
results supporting this conjecture. Most importantly, it is true for
toric log Fano varieties with empty $B$ by Borisovs
\cite{BorisovBorisov92}.
Fano 3-folds with terminal
singularities (i.e.\  1-log terminal) are bounded (see Kawamata
\cite{Kawamata89} for $\mathbb Q$-factorial case, general case is due to
Mori, unpublished). Log Fano 3-folds of fixed index $N$
and Picard number 1
(in particular, they are $(1/N)$-log terminal) are bounded by
\cite{BorisovA93}.
  \end{say}

  \begin{say}
In the case of positive $K$ nothing is known in dimension greater than
two.
  \end{say}

  \begin{say}
Let me list the main points where dimension two is essential in the
proofs: diagram method, lemma~\ref{only surfaces}, lemma~\ref{can
decrease all by delta}. It could also be argued that all proofs
involving $D^{2}$ do not have a chance to be generalized to higher
dimensions because only in dimension $n=2$ $D^{n}$ behaves well under
birational transformations. However, this obstacle can be avoided by
using systematically $\underset{c}\le$ instead of $\le$ and then
applying~\ref{easy inequalities}.
  \end{say}

  \begin{say}
Theorem~\ref{bound for nonsingular when -K nef} has applications to
singularities of curves on surfaces with positive $-K$ such as the
projective plane or Del Pezzo   surfaces.
For example, it implies that an irreducible curve of
degree $d$ for $d\gg0$ cannot have tangents of high multiplicity as
compared to $d$. It would be interesting to spell this connection out
and compare with known results.
  \end{say}

  \begin{say}
Most of the results of this paper remain valid if one works in the
category of algebraic spaces of dimension two
instead  of the category of algebraic surfaces.
  \end{say}

\ifx\undefined\bysame
\newcommand{\bysame}{\leavevmode\hbox to3em{\hrulefill}\,}
\fi

\end{document}